\documentclass[12pt]{article}
\usepackage{amsmath,amssymb,amsthm,euscript} 
\usepackage[cal=esstix,frak=stixtwo,scr=boondox]{mathalfa}
\usepackage{caption}
\usepackage{xcolor,graphicx,epsfig,subfigure}
\usepackage{tikz}
\usepackage[colorlinks = true, urlcolor = blue, pdfborder= 0 0 0.5, urlbordercolor = blue ]{hyperref}
\DeclareFontEncoding{FML}{}{}
\DeclareFontSubstitution{FML}{futm}{m}{it}
\DeclareSymbolFont{myletters}{FML}{futm}{m}{n}
\DeclareMathSymbol{\mapg}{\mathord}{myletters}{000103}

\newsavebox{\foobox}
\newcommand{\slantbox}[2][0]{\mbox{%
    \sbox{\foobox}{#2}%
    \hskip\wd\foobox
    \pdfsave
    \pdfsetmatrix{1 0 #1 1}%
    \llap{\usebox{\foobox}}%
    \pdfrestore
}}
\newcommand\unslant[2][-.25]{\slantbox[#1]{$#2$}}

\newcommand{\uppartial}{\unslant\partial}
\newcommand{\upPartial}{\unslant\partial\kern-0.8pt}

\newtheorem{theorem}{Theorem}
\newtheorem{lemma}{Lemma}
\newtheorem{example}{Example}
\newtheorem{cor}{Corollary}
\newtheorem{conjecture}{Conjecture}
\newtheorem*{notation}{Notation}
\newtheorem{proposition}{Proposition}

\theoremstyle{definition}
\newtheorem{definition}{Definition}
\newtheorem{problem}{Question}
\newtheorem{remark}{Remark}
\renewenvironment{proof}{
\par\noindent{\it Proof.}} { \mbox{}\hfill
$\blacksquare$ \par }

\normallineskip=1pt \normalbaselines \setbox\strutbox\hbox{\vrule
width 0pt height 6.4pt depth 1.6pt}

\makeatletter
\renewcommand{\@oddhead}
{\raisebox{0pt}[\headheight][0pt]{%
\vbox{\hbox to \textwidth{\strut  Mark Pollicott and Polina Vytnova \hfill
\strut Zeta functions zeros \ }%
\hrule}}}
\renewcommand{\@oddfoot}{\hfill ---~\thepage~--- \hfill}
\makeatother

\newcommand\const{\mathord{\rm const}}
\newcommand\trace{\mathop{\rm trace}}

\newcommand\li{\mathord{\rm li}}

\newcommand\eqdef{\buildrel{\rm def}\over=}

\newcommand\Mat{\mathord{\rm Mat}}
\newcommand\Isom{\mathord{\rm Isom}}

\renewcommand\mod{\allowbreak\space\mathord{\rm mod}\allowbreak\, }

\newcommand\tr{\mathop{\rm tr}\nolimits }

\newcommand{\bbC}{\mathbb{C}}

\newcommand\acosh{\mathop{\rm ArcCosh}}
\newcommand\Arg{\mathop{\rm arg}\nolimits }

\mathchardef\R"71B3 \mathchardef\Z"71B4 \mathchardef\C"71B2
\mathchardef\Q"71B1 \mathchardef\N"71B0 \mathchardef\k"71B9

\newcommand\dfill{\cleaders\hbox to 10pt{\hss.\hss}\hfill }

\newcommand\LR[3]{\setbox0\hbox{$#3$}\setbox1\hbox{$\left#1\vcenter{\copy0}
\right#2$}\dimen200\ht1\advance\dimen200 by -\ht0
\dimen201\ht1\advance\dimen201 by \dp1
\mathord{\mathopen{\lower\dimen200\hbox{$\left#1\vcenter to
\dimen201{}\right.$}}\copy0\mathclose{\lower\dimen200
\hbox{$\left#2\vcenter to \dimen201{}\right.$}}}}

\makeatletter
\def\fdsy@scale{1}
\newcommand\fdsy@mweight@normal{Book}
\newcommand\fdsy@mweight@small{Book}
\newcommand\fdsy@bweight@normal{Medium}
\newcommand\fdsy@bweight@small{Medium}
\DeclareFontFamily{U}{FdSymbolA}{}
\DeclareFontShape{U}{FdSymbolA}{m}{n}{
    <-7.1> s * [\fdsy@scale] FdSymbolA-\fdsy@mweight@small
    <7.1-> s * [\fdsy@scale] FdSymbolA-\fdsy@mweight@normal
}{}
\DeclareFontShape{U}{FdSymbolA}{b}{n}{
    <-7.1> s * [\fdsy@scale] FdSymbolA-\fdsy@bweight@small
    <7.1-> s * [\fdsy@scale] FdSymbolA-\fdsy@bweight@normal
}{}
\makeatother

\DeclareSymbolFont{mysymbols}{U}{FdSymbolA}{m}{n}
\SetSymbolFont{mysymbols}{bold}{U}{FdSymbolA}{m}{n}

\DeclareMathSymbol{\medtriangleright}{\mathalpha}{mysymbols}{86}
\DeclareMathSymbol{\medtriangleleft}{\mathalpha}{mysymbols}{88}
\DeclareMathSymbol{\medtriangleup}{\mathalpha}{mysymbols}{87}
\DeclareMathSymbol{\medtriangledown}{\mathalpha}{mysymbols}{89}

\usepackage[backgroundcolor=white!90!red]{todonotes}

\title{The Bowen--Series coding and zeros of zeta functions }
\author{%
  Mark Pollicott \and Polina Vytnova\thanks{
The first author was partly supported by ERC-Advanced Grant $833802$-Resonances
and EPSRC grants EP/W033917/1 and EP/V053663/1. 
The authors would like to thank the organizers of the Simons Semester program 
``Topological, smooth and holomorphic dynamics, ergodic theory, fractals''
for the financial support and hospitality. 
  The authors are also grateful to the organizers of the program at SUSTech in 
Shenszen (China) and to the head of the mathematics department Jeff Xia 
for their help and hospitality.}
}%
\date{ }

\newcounter{mylist}

\begin{document}

\maketitle

\begin{abstract}
We give a discussion of the classical Bowen--Series coding and, in particular,
its application to the study of zeta functions and their zeros. 
In the case of compact surfaces of constant negative curvature $\kappa = -1$ the
analytic extension of the Selberg zeta function to the entire complex plane is
classical, and can be achieved using the Selberg trace formula.  However, an
alternative dynamical approach is to use the Bowen--Series coding on the
boundary at infinity to construct a piecewise analytic expanding map from which the
extension of the zeta function can be obtained using properties of the
associated transfer operator.  This latter method has the advantage that it also
applies in the case of infinite area surfaces provided they don't have cusps.
For such examples the location of the zeros is somewhat more mysterious.
However, in particularly simple cases there is a striking structure to the
zeros when we take appropriate rescaling.  We will try to give some insight into
this phenomenon. 
\end{abstract}

\tableofcontents

\bigskip

\paragraph{Background.}
Some useful references for the basic material in  these lectures are contained in the following
books and articles.   We will not cover even a fraction of  all of the material
in these sources, but that  does not detract from their interest. 
The selection reflects the theme we are trying to give to the story we are telling.

\begin{enumerate}
  \item Intoductory reading: A lot of the foundational material can be found in  a very
    accessible form in the Proceedings of the Trieste Conference on {\it Ergodic Theory, 
    symbolic dynamics and hyperbolic dynamics} from 1989~\cite{BKS}.
\item Ergodic theory: The LMS lecture notes series volume by  Nicholls~\cite{Nicholls}.
\item Trace formulae: The very readable article of McKean on the use of
    heat kernels~\cite{McKean} and a monograph by Buser~\cite{B92}.
\item Resonance-free regions and density of resonances: two articles of
    Naud~\cite{Naud1} and~\cite{Naud2}.
\item Since the notes were originally written, a book by Borthwick~\cite{B18} and a book by
    Dyatlov and Zworski~\cite{DZ22} appeared. 
\end{enumerate}


\section{Introduction}
Dynamical  systems and ergodic theory are both subjects rich in concrete  examples.
To set the scene for these notes we will consider two basic examples of two somewhat 
different types of dynamical systems.  
The first is an example of a discrete dynamical system (or map) associated to
an~$\mathbb Z_+$-action and corresponds to an expanding map on a disjoint union
of intervals.  This example has the virtue of being both simple and accessible.  
The second will be an example of a continuous dynamical system (or flow) associated 
to an $\mathbb R$-action and corresponds to the classical example of the geodesic 
flow for a surface of constant negative curvature. This example has been of paramount 
importance in the original development of ideas in ergodic theory, particularly in the 1930s. 
 
Both of these examples,  the map and the flow,  have been the subject of extensive research. 
However, as we will see later, these two systems are not as unrelated as they might at first appear.
The basis for this connection is some elegant ideas in the work of Bowen and
Series~\cite{BS79} (and a closely related approach appears in the work of Adler and
Flatto~\cite{AF91}).
Unfortunately, the history of the theory is tinged with sadness.  Robert (Rufus)
 Bowen was a professor at the University of California at Berkeley. During
 the collaboration with Caroline Series he died suddenly, at the age of 31.
 Their joint work was completed by his co-author and published posthumously in 
 a memorial volume of Publications Mathematiques (IHES)  dedicated to Bowen. 
 Caroline served as the President of the London Mathematical Society in 2018--2022.

There have been other ways to associate to a geodesic flow an expanding map.
Perhaps the  best known is that of Ratner, which she subsequently generalised to Anosov
flows. Later the construction was further generalised by Bowen to Axiom A flows (which were
introduced by Smale, his supervisor in Berkeley).  This method used two-dimensional 
Markov sections transferred along the flow direction, with obtained ``parallelepipeds''  
playing the role of flow boxes. The associated Markov map corresponds to the interval map 
induced from the Poincar\'e map between sections by collapsing them along the
stable direction (an observation which, at least, is made explicit in a paper of
Ruelle~\cite{R76}).  However, the resulting interval map is not canonical,  as one
would anticipate from its potentially greater generality.   

On the other hand, the Bowen--Series coding is very closely related to the
generators of the fundamental group of the surface. 
Although the choice of generators is not unique, this leads to a very natural 
approach to understanding orbits.  
Indeed, a part of the appeal of the Bowen--Series coding is its transparency.
The work of Bowen and Series pre-dated later work of Cannon and others on the
automatic structure of more general Gromov hyperbolic groups, which is in a
similar spirit.    

The Bowen--Series coding has proved to be very useful in a number of different applications. 
This is particularly true when an explicit understanding of the  fundamental group plays a role.
The Selberg zeta function takes into account precisely one closed geodesic in each conjugacy
class of the fundamental group.

In Section~\ref{s:locres} we shall consider an application of one of the
simplest cases of such coding, for a three funnelled surface or a ``pair of
pants'' and for a one-holed torus, to the study of resonances of the associated Selberg zeta function.
To complete the introduction we summarize what we want to do, in a nutshell.

{\bf Aim of the note:}
 To explore which results/methods (particularly in the context of zeta functions) carry over from
the case of compact surface $V$ to a surface $V$ having infinite area using the Bowen--Series coding. 

\section{Two types of dynamical system}
We now compare these two types of dynamics.  

  \bigskip
\begin{figure}[h]
  \centerline{\includegraphics{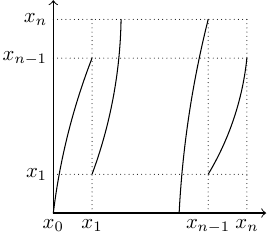}   \rule{30mm}{0pt} 
  \raisebox{24pt}{\includegraphics{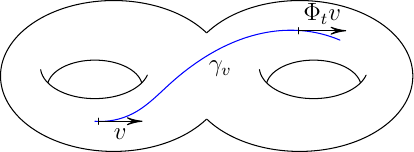}}}
  \centerline{(a) \rule{80mm}{0pt} (b) \rule{10mm}{0pt}}
  \caption{(a)~A Markov expanding map of the interval. (b)~The geodesic flow on the unit 
  tangent bundle of a compact surface of constant negative curvature.}
  \label{mapT:fig}
\end{figure}
  
\bigskip

\begin{table}
\begin{tabular}{||p{81mm}|p{81mm}||}
  \hline
  \multicolumn{1}{||c|}{\bf Discrete (Interval maps)
  \raisebox{-8pt}{\rule{0pt}{24pt}}} &
  \multicolumn{1}{c||}{\bf Continuous (Geodesic flows)} \\ 
  \hline
  Consider a partition of the unit interval $X = \bigcup\limits_{j=0}^{n-1} I_j$, 
  where $I_j=[x_j, x_{j+1}] $ and $0=x_0 < x_1 < \ldots < x_n=1$  & 
  Consider a compact connected surface~$V$ with constant negative curvature $\kappa=-1$ such that
  the first homotopy group~$\pi_1(V)$ is finitely generated. Let 
  $ Y =\{ v \in TV \colon \|v\| = 1 \}$ be the unit tangent bundle. \\
  \hline
 Let  $T \colon X \to X$ be a piecewise smooth, uniformly expanding, Markov and
 transitive map, namely 
   $T \in C^\infty(I_j)$;   
   $\exists \beta>1$ such that $\forall x \in X$
     $|T^\prime(x)|\ge\beta$; $T(I_j)$ is a union of $I_k$; &
 Let $\varphi_t \colon Y \to Y$ be the geodesic flow:
 $\forall v \in Y$ choose the unique geodesic $\gamma_v \colon \mathbb{R} \to V$ such that 
  $\dot\gamma_v(0)= v$ and then define $\varphi_t v = \dot\gamma_v(t)$. 
  This corresponds to parallel transport. \\
   $ x \in X$ which orbit $\{ T^n  x\} \subseteq X$ is dense. &
    The flow~$\varphi_t$ is transitive, i.e. 
  $\exists v \in Y$ which orbit $\{ \varphi_t v \subset Y \mid t \in \mathbb R\} $ is
  dense. \\
  \hline
   {\bf Lemma~1} (Folklore).
      There exists a unique $T$-invariant
    probability measure, which is absolutely continuous with respect to the Lebesgue measure.
  & 
  {\bf Lemma~2.} Let $d\vartheta$ be the Lebesgue measure along the orbits of
  the flow. Then the Liouville measure $d \nu \!=\! d(\mathrm{Vol})_V \times d
  \vartheta$ is invariant with respect to the geodesic flow $\varphi_t$.
  \\
  \hline
\end{tabular}
\caption{Correspondence between discrete and continuous dynamics.}
\label{tab:discrete-continuous}
\end{table}
\bigskip


\setcounter{lemma}{2}

We list corresponding properties of the 1d map and the flow in
Table~\ref{tab:discrete-continuous}.  
The folklore Lemma~1 on existence and uniquence of the invariant measure has
been variously attributed to: Bowen (1979), Adler (1972), Flatto (1969), Weiss
(1968), Sinai (1968) and  Renyi (1960). For an interesting historical
perspective we refer the reader to~\cite{Adler-Bowen}.

Of course there are various classical results we are implicitly using.
  For example, we will assume that the surface $V$ is orientable and has genus
  at least~$2$ (or, equivalently, negative Euler characteristic).  This ensures
  that $V$ can carry a metric of negative curvature.  Moreover, there is a
  number of obvious generalisations that we will not consider.  For
  example, the analogous results in the case of surfaces of variable negative
  curvature. 

Having emphasized the parallels between these two systems, we can look for a
more tangible connection.  More precisely, we can pose the following question. 

\begin{problem}
  How can we relate the Markov map and the geodesic flow?
\end{problem}

One elegant solution to this problem appears in the work of Bowen and
Series~\cite{BS79}, as mentioned in the introduction, where an explicit
construction is suggested. 
  If we take a surface of a fixed genus~$g\geq 2$ then the associated
  expanding maps  for the geodesic flow for different metrics on~$V$ will be
  topologically conjugate.  However, the actual maps themselves will depend
  fundamentally on the choice of metric on the surface.   
 
  Recall that the universal cover of a compact connected surface~$V$ that carries a metric of constant
  negative curvature is the unit disk~$\mathbb D$ endowed with the hyperbolic
  metric; and the fundamental group $\pi_1(V)$ acts on the universal
  cover~$\mathbb D$. The Bowen--Series construction rests on the idea of hyperfineteness,
  i.e., the way the action 
   $$
   \pi_1(V) \times \uppartial \mathbb D \to\uppartial \mathbb D 
   $$
   of the fundamental group~$\pi_1(V)$ on the boundary~$\uppartial\mathbb D$ of
   the universal cover can be replaced by a single transformation.
   The Bowen--Series transformation is  a specific realisation of such a
   phenomenon.
It is defined as a map on the limit set, and more precisely is defined on a union of
arcs which form a partition of the limit set.  However, we can
easily think of them as being maps of the interval, or maps on a Cantor set
contained in an interval.   

This approach has the advantage that complicated problems for the geodesic
flow can often be reduced to simpler problems for the expanding map of the
interval.  We can summarise this method (elevated below to a ``philosophy'') as
follows.

\fbox{\parbox{0.98\textwidth}{
\noindent {\bf Philosophy} 
\begin{enumerate}
        \parsep0pt
    \itemsep0pt
  \item Consider a problem for the geodesic flow $\varphi_t \colon Y \to Y$; 
  \item Reduce the problem to a problem for the associated expanding map  $T \colon X \to X $;
  \item Solve the problem for $T$ (assuming one can!);
  \item Relate it back to a solution for the original problem for $\varphi_t \colon Y \to Y$.
\end{enumerate}
}}

\bigskip 

\noindent Of course, like all approaches its value probably ultimately depends on how
useful it is.  More precisely, we could ask:   

\begin{problem} 
  What sort of problems can one address?
\end{problem}
  
One would naturally expect any useful technique to have many different applications, but for simplicity we could divide 
the types of problems we  typically consider  into two general classes: \\
\noindent{\bf Types of problems} \quad 
\begin{enumerate}
  \item[(a)] {\it Topological}, such as properties of closed orbits).  As a specific
    example we could consider results on the distribution of closed orbits which
    reflect information on their free homotopy classes (i.e., conjugacy classes
    in the fundamental group).  The Bowen--Series coding makes it easy to relate
    the closed orbits for the flow to periodic orbits $\{x, Tx, \cdots, T^{n-1}x
    \}$ (with $T^nx=x$) for the associated transformation.  Furthermore, the
    usefulness of this particular coding is that the word length of the
    corresponding geodesic will be $n$, except in a finite number of exceptional
    cases.
  \item[(b)] {\it Measure theoretic}, such as properties of the invariant
      measures for the map and for the flow. The two measures are intimately 
      related by the classical work of Bowen--Ruelle. For example, showing
      ergodicity for the flow automatically implies ergodicity for the discrete
      map.  In the reverse direction, ergodicity for the discrete map, together
      with some modest hypothesis on the roof function, implies ergodicity for the flow. 
      Stronger properties such as strong mixing, mixing rates, central limit
      theorems, etc. can be considered in each context, with varying degrees of
      complexity for the correspondence between them.  
\end{enumerate}

\noindent Since we are interested in ergodic theory and dynamical systems it is
natural to use the model of the geodesic flow using expanding interval maps to
try to understand its dynamics.  However, in the interests of being open minded,
we should also ask: 

\begin{problem}
  Is this the best approach? 
\end{problem}

The word ``best'' is a little subjective, so there is no real definitive answer.   
However, the somewhat equivocal answer could be ``sometimes, yes'', however, of course, 
this depends on what we are interested in and what we mean by the question. 
For example if we are interested in closed geodesics on the surface then these
can be viewed dynamically as closed orbits for the geodesic flow.  At the same
time, in some cases their study can be advanced by other (less dynamical) methods.
For example, we note that  
\begin{enumerate}
  \item If~$V$ is a compact surface of constant negative curvature there are
      already powerful techniques (e.g. trace formulae, representation theory)
      which can often give more precise results.   
  \item Nevertheless, if $V$ is a non-compact surface of infinite area the
      above  methods may not work so  well, but  the dynamical approach often 
    still applies.  At a certain level this can be thought of as being because
    the dynamical method uses only the compact recurrent set of the geodesic
    flow. 
\end{enumerate}

\noindent Let us now recall a classical example of a compact surface.

\begin{example}
Let~$V$ be an oriented connected compact surface of genus $g \geq 2$.  This surface not
only supports a metric of constant curvature $\kappa=-1$, but the space of such
metrics is  $(6g - 6)$-dimenional.
\end{example}

Next we want to consider a couple of examples of infinite area surfaces with constant curvature $\kappa=-1$ 
to illustrate the second part.

\begin{figure}[h]
  \centerline{
  \includegraphics[height=40mm]{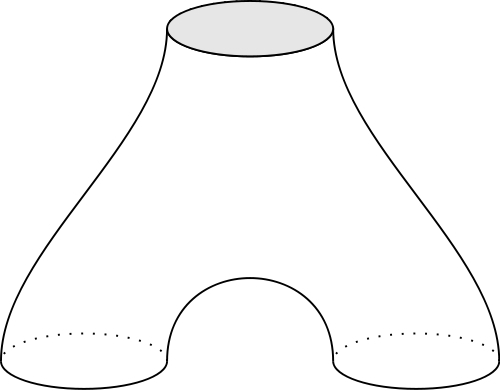}
  \rule{20mm}{0pt}
  \raisebox{5mm}{\includegraphics[height=30mm]{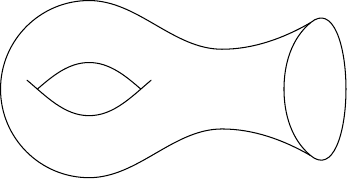}  }
  }
  \caption{Left: A pair of pants, or a 3-funnelled surface. Right: A one-holed torus or a torus with a funnel. 
  Both drawings represent an artistic impression. According to the Hilbert Theorem,
  no complete surface of constant negative Gaussian curvature can be
  isometrically immersed in~$\mathbb R^3$. }
  \label{pants:fig}
  \end{figure}

\begin{example}\label{pants:ex}
 A three funnelled surface, or ``a pair of pants'', is a surface homeomorphic to a sphere with three disjoint
 closed disks removed. It carries metrics with constant negative curvature
  $\kappa=-1$.  We can either consider the complete surface which has three
  funnels (or ends) or alternatively we can cut along three closed
  \emph{shortest} geodesics around each of the funnels to get a surface with three boundary components (as
  in Figure~\ref{pants:fig}, left).  In point of fact, from a dynamical viewpoint the distinction
  is unimportant. For the geodesic flow on either version the important
  dynamical component is the recurrent part of the geodesic flow on the
  $3$-dimensional unit tangent bundle. It is a compact set whose intersection
  with any transverse two-dimensional set is a Cantor set. 
  Since none of the orbits can cross any of the three geodesics around the 
  funnels, the corresponding set appears on the pair of pants. 
 
\end{example}
 
\begin{example}
  \label{torus:ex}
A one-holed torus, or a torus with a funnel, is a surface homeomorphic to a
torus with a closed disk removed. Similarly to a pair of pants, it carries a metric with constant negative
curvature~$\kappa=-1$.  Again, we can either consider the complete surface which has a
single funnel (or end) or alternatively we can cut along the closed geodesic
around the funnel to get a surface with a single  boundary component (as in
Figure~\ref{pants:fig}, right). As in Example~\ref{pants:ex}, from a dynamical viewpoint the
distinction is unimportant.  The important dynamical component is the recurrent
part of the geodesic flow and  none of the orbits in the recurrent part can
cross the  geodesic around the funnel. 
\end{example}

We would like to note that although the approach of Bowen and Series is very
successful in the case of surfaces, it does not naturally generalise to higher dimensions.

\section{Hyperbolic Geometry} 
In this section we want to introduce the setting in which we will be working.
This involves recalling  some basic definitions and notions in hyperbolic geometry.
Hyperbolic geometry was developed  in the 19th century by Gauss and
Bolyai\footnote{This geometry was independently developed by a Russian
mathematician N.~I. Lobachevski, and carries his name in post-Soviet
countries.} to show that the 5th postulate of Euclid was not implied by the others. 

There are many different models for the hyperbolic space, but we will concentrate on the disk model.  
We begin with the basic definitions, and refer the reader to the book of
Beardon~\cite{B83} for more details.

Let us denote by $  \mathbb D = \{ z \in \mathbb C \mid |z|<1 \} $
   the open unit disk in the complex plane.  
\begin{definition}
   We can equip~$\mathbb D$  with the 
  {\it Poincar\'e metric} defined locally by 
  \begin{equation}
    ds = \frac{2 |dz|}{1 - |z|^2}.
  \label{eq:hypmetric}
  \end{equation}
\end{definition}
\noindent In particular, the distance between two points $z_1, z_2 \in \mathbb
D$ is 
$$
d(z_1, z_2) = 2 \tanh^{-1} \left| \frac{z_1 - z_2}{1 - z_1 \overline z_2} \right|.
$$
This defines a complete Riemannian metric on~$\mathbb D$. The factor of~$2$ in
the definition of the metric~\eqref{eq:hypmetric} is ensuring a convenient normalisation for the curvature.
The Poincar\'e metric has several useful properties which we briefly summarize
below.
\begin{lemma}
  \label{pmetric1:lem}
Consider the unit disk equipped with the Poincar\'e metric $(\mathbb D,ds)$. Then:
\begin{enumerate}
  \item The space has Gaussian curvature $\kappa=-1$; and 
  \item The orientation preserving isometries take the form of linear
    fractional transformations 
    $$
      \mapg \colon \mathbb D \to \mathbb D; \quad  
      \mapg(z) = \cfrac{\alpha_1 z+\alpha_2}{\overline{\alpha_2} z + 
      \overline{\alpha_1}}, \quad  \mbox{ where } |\alpha_1|^2 - |\alpha_2|^2=1, \,
      \alpha_1, \alpha_2 \in \mathbb C .
    $$
\end{enumerate}
\end{lemma}
 \noindent  Note an ambiguity with respect to the simultaneous sign change
  $(\alpha_1,\alpha_2)
    \longleftrightarrow (-\alpha_1, -\alpha_2)$ which one has to bear in mind.
There are other equivalent models for the hyperbolic plane, such as the Poincar\'e
upper half plane model.   However the symmetry of the disk model is particularly
useful in what follows in later sections.

For a compact surface~$V$ with negative curvature $\kappa = -1$, the
Poincar\'e metric on the unit disk has a particular significance.

\begin{lemma}
\label{pmetric2:lem}
Let~$V$ be a compact surface of constant  curvature $\kappa=-1$. Then:\\
\begin{minipage}{0.725\textwidth}
\begin{enumerate}     
  \item The Gauss--Bonnet theorem implies that the genus~$g$
      of~$V$ satisfies $g \ge 2$. Equivalently, the Euler characteristic is
      strictly negative.   
  \item The inequality $g \geq 2$ gives that the universal cover of~$V$
      is~$\mathbb D$ and the lifted metric is the Poincar\'e metric
      cf.~\cite{B83}. 
  \item The geodesics on~$V$ (locally distance minimizing) lift to geodesics
      on~$\mathbb D$. The latter are circular arcs which meet the unit
      circle~$\uppartial\mathbb D$ orthogonally. 
\end{enumerate}
\end{minipage}
\begin{minipage}{0.02\textwidth} \quad
\end{minipage}
\begin{minipage}{0.25\textwidth}
\includegraphics{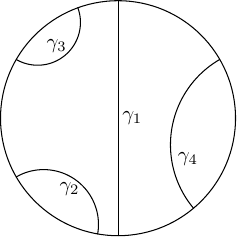}

Figure 3: The geodesics in the disk model of the hyperbolic plane. 
\end{minipage}
\end{lemma}

We have identified in Lemma~\ref{pmetric1:lem} the group of isometries for the Poincar\'e metric.  
We next want to consider its discrete subgroups. 

\section{Fuchsian groups and the Bowen--Series map}
The next ingredient to consider is a suitable subgroup of the isometries of the
unit disk with respect to the Poincar\'e metric. Given a compact surface
of constant curvature $\kappa = -1$, we  can write $V = \mathbb D/\Gamma$ where 
$\Gamma < {\rm Isom} (\mathbb D,ds)$ is a discrete subgroup  isomorphic to
the fundamental group of~$V$, that is $\Gamma =  \pi_1(V)$.

\subsection{Fuchsian groups}
We recall the following standard definition of a Fuchsian group.

\begin{definition}
  We say that $\Gamma < \Isom (\mathbb D,ds)$ is a {\it Fuchsian group} if it is a
  discrete subgroup. 
\end{definition}
\noindent  In other words,~$\Gamma< \Isom (\mathbb D,ds)$ is Fuchsian, if the orbit $\Gamma 0$ of the
point $0 \in \mathbb D$ under the action of~$\Gamma$ doesn't have accumulation points inside the disk. 
Moreover, since $\mapg \in \Gamma$ is a
linear fractional transformation, it can be extended to the boundary $\uppartial
\mathbb D$, furthermore, $\mapg$ preserves the boundary. 
Hence we can associate to the Fuchsian group $\Gamma$ an important closed subset
of the boundary circle $\uppartial\mathbb D$ called the limit set.

\begin{definition}
The limit set $\Lambda \subset \uppartial\mathbb D$ consists of Euclidean accumulation
points of the orbit of the centre~$0$, namely,
$$
\Lambda = \overline {\Gamma 0} - \Gamma  0.
$$ 
\end{definition}
\noindent It is not essential to take the orbit of $0$; any point in the interior of
$\mathbb D$ would do.
The limit set is always compact.   Moreover we know the following:

\begin{lemma}
  \label{limdim:lem}
  The limit set of a Fuchsian group~$\Gamma$ is either:
\begin{enumerate}
\item The entire boundary circle, $ \Lambda = \uppartial\mathbb D$; or
\item A Cantor set of Hausdorff dimension $0 < \delta < 1$.
\end{enumerate}
\end{lemma}
\noindent These are usually referred to as limit sets of Type~I and Type~II, respectively.

For the moment, we will focus on the case when~$V =\mathbb D/\Gamma $ is a compact surface and the associated 
Fuchsian group~$\Gamma$ is called co-compact. In this case the limit set is the entire circle.
The elements of the Fuchsian group~$\Gamma$ have a few useful and important properties.

\begin{lemma}
Let~$V = \mathbb D/\Gamma$ be a compact connected surface of constant negative
curvature and genus~$g$. 
\begin{enumerate}
  \item  Every transformation $\mapg \in \Gamma\setminus\{e\}$ has two
    fixed points, a repelling and an attracting one, on $\uppartial\mathbb D$.
    We denote the {\it attracting point} by $p_{\mapg }^+$, so that
    $|{\mapg }^\prime (p_{\mapg }^+)|<1$.  Similarly, we
    let  $p_{\mapg }^-$ denote the {\it repelling point}, so that $|{\mapg }^\prime(p_{\mapg }^-)|>1$.  
  
  \item The group $\Gamma$ is finitely generated and finitely presented. 
      For instance, there exists a finite set of generators $a_j \in \Gamma$, 
      $j = 1, \ldots, 2g$ such that  
    $$
    \Gamma = \langle a_1,a_2,\ldots, a_{2g} \mid (a_1 a_2^{-1} a_3
    a_4^{-1} \ldots a_{2g-1} a_{2g}^{-1} ) \cdot (a_1^{-1} a_2 a_3^{-1} a_4 \ldots
    a_{2g-1}^{-1} a_{2g}) = 1 \rangle.
    $$
\end{enumerate}
\end{lemma}

\subsection{Bowen--Series map}
Given a Fuchsian group we want to define a piecewise continuous expanding Markov map 
called the {\it Bowen--Series map} defined on the boundary $\uppartial\mathbb D$.  

   In order to define the map $T \colon \uppartial\mathbb D \to
   \uppartial\mathbb D$ we need to define a partition of the boundary. The
   construction below relies on Ford fundamental domain.    
   
   To every $\mapg \in \Gamma$ we relate an \emph{isometric circle}:
    $$  
    C(\mapg ) = \{z \in \mathbb D \mid |\mapg ^\prime(z)|=1 \} 
    $$
  The name isometric here refers to the fact that $\mapg$ locally doesn't change
  Euclidean length. A nice treatment of isometric circles can be found
  in Ford's monograph~\cite{Ford19}.

  The isometric circle of an isometry $C(\mapg )$ is a geodesic 
  in the Poincar\'e metric. Indeed, assuming $\mapg (z) =
  \frac{\alpha_1 z+\alpha_2}{\overline{\alpha_2} z + \overline{\alpha_1}} $ we
  get $\mapg^\prime(z) = (\overline{\alpha_2} z +\overline{\alpha_1})^{-2}$.
  Therefore $|\mapg^\prime(z)|=1$ defines a circle of radius
  $|\overline{\alpha_2}|^{-1}$ centred at $-\overline{\alpha_1}/\overline{\alpha_2}$ provided $\alpha_2 \ne 0$. 
  To see that it is a geodesic, observe that the triangle with vertices at the
  origin, the centre of $C(\mapg)$, and a point of the intersection of $C(\mapg)$ and
  $\uppartial \mathbb D$ is right-angled. 
  A useful property is that a map $\mapg \in \Gamma$ transfers its isometric circle to the isometric
  circle of its inverse: $\mapg  C(\mapg ) = C(\mapg ^{-1})$.  

  The region of the
  unit disk ``exterior'' to the isometric circles of \emph{all} elements
  of~$\Gamma$ is a fundamental domain for the Fuchsian group~$\Gamma$. It is
  called~\emph{Ford domain}. It is known that every side $s$ of the fundamental
  domain is identified with another one by an element~of~$\Gamma$ that we shall
  denote~$\mapg(s)$. Moreover,
  $$
  \Gamma = \langle \mapg(s) \mid s \mbox{ is a side of the Ford fundamental domain of } \Gamma
  \rangle,
  $$
  in other words, the elements identifying the sides of the Ford domain
  constitute a set of generators of~$\Gamma$. 
  
  In order to guarantee that the boundary map we are about to define is Markov, we need
  to impose additional assumptions on the way the surface group acts on 
  the hyperbolic plane:

  \begin{itemize}
      \item Every side~$s$ of a fundamental domain lies on the isometric circle
          of~$\mapg(s)$, i.e. $s \subset C(\mapg(s))$;
      \item All hyperbolic lines, containing sides of the fundamental domain,
          are entirely covered by
          images of the sides of the fundamental domain under action
          of~$\Gamma$.   
  \end{itemize}

\begin{remark}
In many cases the assumption may not hold, but then there exists another
representation of the fundamental group $\Gamma=\pi_1(V)$ for which this property holds.  
The choice of~$\Gamma$ was written down independently by both
Bowen--Series~\cite{BS79} and Adler--Flatto~\cite{AF91}.
In the interests of clarity and historical reverence we will use the
original construction of Bowen and Series~\cite{BS79}.
\end{remark}

\noindent\begin{minipage}{0.68\textwidth}
Under the above assumptions, configuration of the isometric circles takes the form illustrated in 
Figure~\ref{fund:fig}, adapted from~\cite[Fig.~3]{BS79}.
The boundary of the fundamental domain is given by arcs of isometric circles of
generators of~$\Gamma$. Each of the circles may intersect with its direct neighbours, but no more.
Their endpoints give a partition of the boundary. 
We can now pick every second point and use these points on the boundary 
circle to divide it into the same number of arcs, and then define an analytic map 
to each of these arcs.
\end{minipage}
\begin{minipage}{0.3\textwidth}
   \centerline{  \includegraphics{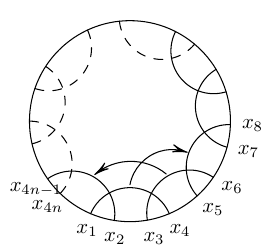} }
   \captionof{figure}{Isometric circles.}
\label{fund:fig}
\end{minipage}

\bigskip

\noindent  This will give us a required map of the circle, that we
treat as an interval with endpoints identified, and the arcs are thought of
as intervals, in the usual way. 
More precisely, we enumerate the sides of the fundamental domain as
$s_1,\ldots,s_{2n}$ in anti-clockwise direction and denote the end points of the isometric circle
corresponding to $s_j$ by $x_{2j-1}$ and $x_{2j+2}$, for $j=1,\ldots,2n-1$, and the end points 
of $s_{2n}$ we label $x_{4n-1}$ and $x_2$. Then we define $T \colon \uppartial
\mathbb D \to \uppartial\mathbb D$ by 
\begin{equation}
T(x) =  \begin{cases} 
        \mapg(s_j)(x), &\mbox{if } x \in [x_{2j-1}, x_{2j+1}), \mbox{ for some } j
        < 2n ; \\ 
        \mapg(s_{2n})(x), &\mbox{if } x \in [x_{4n-1}, x_1).
    \end{cases}
  \label{mapTdef:eq}
\end{equation}   
    
We summarise the main properties of the map~$T$ in the following lemma.
\begin{lemma}
 The interval map~$T$ is: 
\begin{enumerate}
\item Real analytic on each interval of the partition specified
    in~\eqref{mapTdef:eq}; 
\item Uniformly expanding: $\exists \beta >1$ such that $|T^\prime(z)|\ge
    \beta$ for every $z \in  \uppartial\mathbb D $;
\item Markov: there exists a \emph{finite refinement} of the partition described
    in~\eqref{mapTdef:eq} with respect to which~$T$ has Markov property. More
    precisely, one needs to take the endpoints of \emph{all} sides of the tesselation
    that go through the vertices of the fundamental domain. The assumptions on
    the action of~$\Gamma$ guarantee the finiteness of the refined partition.
\end{enumerate}
\end{lemma}
\noindent These properties are very classical and lead to a clear understanding of the dynamics of the map.
Moreover, these can be used to understand the induced action of the Fuchsian group~$\Gamma$ on the boundary.
In particular, we  would like to replace the action 
$$
\Gamma \times \uppartial\mathbb D \to \uppartial\mathbb D \qquad (\mapg,x) \mapsto
\mapg(x)
$$
by a single transformation $T \colon \uppartial\mathbb D \to \uppartial\mathbb
D$ (up to finitely many points, which are the end points of the isometric
circles). 

Let us denote by $\Gamma x$ the orbit of $x$ under the action of the group, i.e., 
$\Gamma x = \{ \mapg x \mid \mapg  \in \Gamma \}$. 
\begin{lemma}
    The orbit~$\Gamma x$ of a \emph{typical} point~$x \in \uppartial \mathbb D$ can be written as an equivalence class 
$$
\Gamma x = \{ y \mid \exists n,m \ge 0 \colon T^n x = T^m y \},
$$
where $T$ is defined by~\eqref{mapTdef:eq}.
\end{lemma}

\noindent\begin{minipage}{95mm}
\begin{example}
Consider a regular hyperbolic octagon with angles~$\frac\pi4$. We can enumerate 
the edges $s_j$, $j = 1, \ldots, 8$. Then we may consider the orientation
preserving isometries~$\Gamma_0$ of the hyperbolic plane which identify 
the opposite sides. Let~$\Gamma$ be the group generated by~$\Gamma_0$. One can
check that the octagon is the fundamental domain of~$\Gamma$ and 
$s_i$ is isometric circle of $\mapg(s_i)$. 
It follows by the symmetry that the Bowen--Series map~$T$ has~$16$ arcs of
continuity, all of the same length. We would like to stress that~$T$ is not
Markov with respect to this partition, one needs to refine it.  (See figure on the right, corresponding to
the case when $V = \mathbb
D/\Gamma$ is the Bolza surface.)
\end{example}
\end{minipage}%
\begin{minipage}{4mm}
\ 
\end{minipage}%
\begin{minipage}{70mm}
    \includegraphics{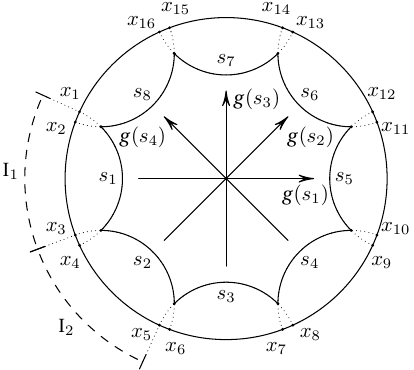}
\end{minipage}

\section{Dynamical Applications:  Compact surfaces}
We now return to the study of the dynamical properties of the geodesic flows on
surfaces, particularly those which can be approached from the perspective of the associated
Bowen--Series map on the unit circle (or more accurately, on a disjoint union
of intervals).  We would like to highlight two particularly well known properties.

\subsection{Counting closed geodesics}
We begin with a topological result. Let~$V$ be a compact surface of constant
negative curvature.   It is easy to see that there are infinitely many closed geodesics
since each  free homotopy class (i.e., conjugacy class) in $\pi_1(V)$ contains
a closed geodesic.  Furthermore, there are countably many closed geodesics since
each free homotopy class for a negatively curved surface contains precisely one
closed geodesic. By a closed geodesic we mean the directed closed geodesics, i.e.,
each curve considered as a set actually corresponds to two directed geodesics,
which differ by orientation.

Let us denote by $\ell(\gamma)$ the length a (primitive) closed
geodesic~$\gamma$. Closed geodesics are in bijection with the (prime)
closed orbits of the geodesics flow of the same length.
   
\begin{notation}
 The number of closed geodesics whose length is at most~$t$  
we denote by
 $$
 N(t) : = \#\{\gamma \mbox{ is primitive closed geodesic and }  \ell(\gamma) \le t\}
 $$
\end{notation}
\noindent It is well known that~$N(t)$ is finite~\cite[Ch 9, \S2]{B92}.
 Let us recall the further properties.
\begin{lemma} 
\label{counting}
  For a compact negatively curved surface we have:
  \begin{enumerate}
     \item $N(t)$ is monotone increasing; and
     \item $\lim\limits_{t\to \infty} \frac1t \log N(t) = 1$.
  \end{enumerate}
\end{lemma}
\noindent The quantity that appears in the second part is the topological
entropy of the geodesic flow. A much stronger result than Part~(2) of 
Lemma~\ref{counting} is the following asymptotic.
\begin{theorem}[Huber~{\cite[Satz 9]{H59}}] Let~$V$ be a compact surface of
    constant negative curvature~$-1$. Then
  \label{huber}
  $$
  \lim\limits_{t \to   \infty}\frac{N(t)}{e^t} \cdot t =1.
  $$
\end{theorem}
      
The natural proof of Huber's theorem  uses the Selberg zeta function $Z_V(s)$
and the location of its zeros. We recall the definition.       
       
\begin{definition}
 We formally define the Selberg zeta function as an infinite product
 \begin{equation}
 Z_V(s) = \prod_{n=0}^\infty \prod_\gamma \left( 1- e^{-(s+n)\ell(\gamma)} \right), \quad s \in \mathbb C,
 \label{zeta:def}
 \end{equation}
 where the product is taken over all closed geodesics on~$V$ and
 $\ell(\gamma)$ stands for the length of~$\gamma$. The
 function is well defined, as the product converges to an analytic function for $\Re(s)>1$. 
\end{definition}
            
The convergence  of the infinite product from~\eqref{zeta:def} on $\Re(s) > 1$
follows from the estimate in Lemma~\ref{counting}. The asymptotic estimates for
the counting function~$N(t)$ follow completely by analogy with the Prime Number
Theorem for primes and the zeta function~$Z_V(s)$ is used in place of the 
Riemann~$\zeta$-function. 

The zeta function $Z_V(s)$ formally defined on~$\Re(s)>1$ can be extended to the
entire complex plane using the famous Selberg trace formula~\cite{McKean}. 
This is explained in the book of Hejhal~\cite{hejhal} in considerable detail.
We recall the basic properties:
\begin{lemma}
 Let $V$ be a compact surface of constant negative curvature. Then the Selberg
 zeta function has the following properties: 
 \begin{enumerate}
    \item $Z_V(s)$ analytic and non-zero for~$\Re(s)>1$; 
    \item $Z_V(s)$ has a simple zero  at~$s=1$;
    \item $Z_V(s)$ has an analytic extension to~$\mathbb C$.
  \end{enumerate}
\end{lemma}
   
\begin{remark}
    Margulis~\cite{M04} gave a generalisation of Huber's theorem to the case of surfaces of
variable negative curvature in his 1972 PhD thesis, that was published in
English more than 30 years later. Margulis' approach to Huber's theorem (for $\kappa
< 0$) uses strong mixing of a suitable invariant measure (nowadays called the
Bowen--Margulis measure) which maximises the entropy. The approach introduced by
Margulis has now been adapted to a number of interesting generalizations. 
\end{remark}   
       
We want to consider error terms in the asymptotic asymptotic relation from 
Theorem~\ref{huber}. To put this into perspective we  need to recall what might
happen in the corresponding setting in number theory. 
       
       
\subsection{Taking a leaf out of the number theory book}
        
Let us begin by recalling some classical ideas from number theory.  These are
useful in understanding the role of the dynamical analogues.

\bigskip
\noindent
{\bf Riemann $\zeta$-function. }
We recall the definition of the Riemann $\zeta$-function: 
    $$
    \zeta(s) = \sum_{n=1}^\infty \frac{1}{n^s} = \prod_p (1-p^{-s})^{-1} 
    $$
where the Euler product is over the prime numbers. Below are some of the
properties of this famous function. 

\begin{lemma}
The Riemann zeta~$\zeta$ function has the following properties:
\begin{enumerate}
  \item $\zeta(s)$ is analytic and non-zero for $\Re(s)>1$; 
  \item $\zeta(s)$ has a simple pole at $s=1$; 
  \item $\zeta(s)$ has analytic extension to $\mathbb C\setminus\{1\}$.
\end{enumerate}
\end{lemma}

Of course there is additionally another property which has not yet
been established.   

\paragraph{Riemann hypothesis (1859):} All zeros of~$\zeta$ in semi-plane $\Re(s)>0$
lie on a line $\Re(s)=\frac12$. 

\medskip

This question is also known as Hilbert's 8'th Problem and a
Clay Institute Millennium Problem. 
There have been many attempts to prove this result, but so far the conjecture
has resisted attempts.  One particularly popular approach is the following:
\paragraph{Hilbert--Polya approach:} Try to relate the zeros of~$\zeta$ to
eigenvalues of some self-adjoint operator.

\medskip

The motivation for this idea is the fact that the spectrum of a self-adjoint operator is
necessarily contained in~$\mathbb R$. Whereas it has proved ellusive for the
Riemann zeta function, the analogue for the Selberg zeta function has been more
successful. 

\subsection{Error terms in counting closed geodesics} 
There is an improvement to Huber's basic  Theorem~\ref{huber} at
least as stated. In fact, the original result of  Huber included error terms, so
now we are actually stating more of his  original theorem.

\begin{theorem}[Huber~{\cite[Satz IV]{H61}}] 
  \label{counterror}
  Let~$V$ be a compact surface with constant gaussian curvature $\kappa = -1$. There exists $\varepsilon>0$ such that 
      $$
      N(t) = \int_2^{e^t} \frac{1}{\log u} {\mathrm d} u (1+O(e^{-\varepsilon t})) =
       \li(e^t) (1+ o (e^{-\varepsilon t})) \quad \mbox{ as } t \to
      +\infty.
      $$
      In other words, there exist $\varepsilon>0$ and $C>0$ such that $\left|\frac{N(t)}{\li(e^t)} - 1 \right| \le C
      e^{-\varepsilon t} $.     
    \end{theorem}
Note that $\varepsilon$ here is responsible for the rate of convergence, rather
then ``just'' for the error term itself.

As is well known, 
$$
\li(x) \colon = \int_2^x \frac{ {\mathrm d} u}{\log u} \sim \frac{x}{\log x} \mbox{ as }x \to +\infty.
$$
Therefore this statement is consistent with the original statement of
Theorem~\ref{huber},  but now  has an additional  error term. 

 The proof uses properties of~$Z_V$  and the rate $\varepsilon>0$ depends on location of
 zeros. This is completely analogous to the  situation in number theory where
 one studies a counting problem for prime numbers by using the Riemann zeta
 function, except in the present context we have stronger results on the zeta
 function. 
        
For compact surfaces there is a natural self-adjoint operator. The extension and
zeros of Selberg zeta function are related to {\it the Laplacian} $\Delta \colon
\mathcal L^2(V) \to \mathcal L^2(V)$.  
This can be defined  as an operator on real analytic functions $\Delta \colon
C^\omega(V) \to C^\omega(V)$ and then extended to  square integrable functions.  
Let 
$$
\Delta \psi_n + \lambda_n \psi_n =0
$$ 
be the eigenvalue equation. There are infinitely many eigenvalues 
$$
  0 = \lambda_0 < \lambda_1\le \lambda_2 \le \ldots \nearrow + \infty
$$
for the operator $-\Delta$.
Moreover, we have the following results:  
  \begin{enumerate}
    \item $\# \{\lambda_n \le t\} \sim \frac{ \mathrm{ Area} (V) }{4 \pi} \cdot t $ as $ t \to \infty$ (Weyl, 1911); 
    \item The zeros of $Z_V(s)$ in $0< \Re(s) < 1$ satisfy the
        equality $s_n(1-s_n) = \lambda_n$, where~$\lambda_n$ are eigenvalues of the Laplacian~\cite{S56}. 
  \end{enumerate}

\noindent\begin{minipage}{113mm}  
The second statement implies that $s_n = \frac12 \pm \sqrt{\frac14 - \lambda_n}$, 
and lie on $[0,1] \cup \{\Re(s) = \frac12 \}$. 
 This can be thought of as an analogue of the ``Riemann Hypothesis''.  
We can formulate it as follows:
\begin{cor} Let~$V$ be a compact surface and let~$s_n$ be a zero of the
    associated Selberg zeta function $Z_V$~defined by~\eqref{zeta:def}.
  Then either $s_n \in [0,1]$ or $\Re(s_n)=\frac12$. 
\end{cor}
 The proof uses trace formulae to relate the zeros of~$Z_V$ to the
 eigenvalues of the Laplacian. It is the self-adjointness of the Laplacian
 which ultimately leads to this property.

\end{minipage}%
\begin{minipage}{5mm}   \quad                  
\end{minipage}%
\begin{minipage}{50mm}                     
    \includegraphics{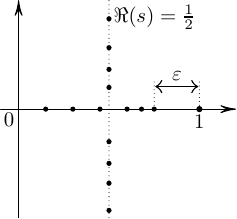}
    \captionof{figure}{The zeros of the Selberg zeta function}
\end{minipage}    

\bigskip                     

The parameter $\varepsilon  > 0$ for the error term for counting closed orbits in
Theorem~\ref{counterror} can be expressed in terms of
the eigenvalues of the Laplacian, namely, it corresponds
to the distance of the zeros of~$Z_V$ from the line~$\Re(s) = 1$ (although
one may need to take~$\varepsilon$ slightly smaller). 

\begin{remark}
The value of $\varepsilon>0$ can be arbitrary small.  This follows from its
spectral interpretation (in terms of the Laplacian) and corresponds to having
eigenvalues close to~$0$. Using classical results of Cheeger~\cite{C70},
Schoen--Wolpert--Yau~\cite{SWY80}
one can show that $\varepsilon \asymp \ell(\gamma_0)$ where $\ell(\gamma_0)$ is
the length of the shortest closed dividing geodesic on~$V$.  
\end{remark}

    The approach and the statement fail if~$V$ has infinite
area. In general case there is no reason for zeros to lie on such lines. 
In the context of non-compact surfaces we need an alternative approach using the Bowen--Series
map, that allows us to establish a correspondence between closed geodesics and closed orbits. 

Recall that a geodesic is called prime if it traces its path exactly once. We
say that a periodic orbit $ \{x,Tx,\ldots,T^{n-1} x \mid T^n x=x \}$ is prime if 
$T^k x \neq x$ for all $1 \le k < n$.

\begin{lemma}\label{closed}  
  There is a bijection between
     closed (prime) geodesics~$\gamma$ on~$V$ and 
     closed (prime) orbits $ \{x,Tx,\ldots,T^{n-1} x \mid T^n x=x \}$,
   where~$T$ is defined by~\eqref{mapTdef:eq}, 
except for a finite number of geodesics corresponding to the endpoints of the
intervals of analyticity of~$T$.
\end{lemma}
            
An advantage of the Bowen--Series construction is that the map~$T$ defined
by~\eqref{mapTdef:eq} can be used to recover the lengths of closed geodesics.  
More precisely,            
\begin{lemma}\label{closed1}
   If~$\gamma$ is a prime closed  geodesic of length~$\ell(\gamma)$ with associated
   prime periodic $T$-orbit 
   $$ \left\{x,Tx,\ldots,T^{n-1} x \mid T^n x=x \right\} ,  
   $$
   then $\ell(\gamma) = \log|(T^n)^\prime(x) |$. 
\end{lemma}
\noindent 
In particular, we can interpret the length of~$\gamma$ as the Lyapunov exponent
of the periodic orbit, which gives a clue to the Lemma. 

To take advantage of this correspondence to study zeta function $Z_V$, we
introduce the family of complex Ruelle--Perron--Frobenius transfer operators
(indexed by $s\in\mathbb C$).  This is often referred to as the ``Thermodynamic
Viewpoint''~\cite{PP90}.   We begin by recalling the definition.

\begin{definition}
Given a (possibly infinite) disjoint union of intervals~$X$ and associated piece-wise
analytic map~$T$ we can define a family of {\it transfer operators} on a space of continuous functions
$\mathcal L_s \colon C^0(X) \to C^0(X)$  by
$$
[\mathcal L_s f](x) = \sum_{Ty=x} e^{-s \log|T^\prime(y)|} f(y), \quad s\in \mathbb C
$$
\end{definition}

\begin{remark}
When $s=1$, this reduces to the ``usual'' Ruelle--Perron--Frobenius operator
for an expanding interval map. 
\end{remark}

However, in order to proceed we need to replace the big space~$C^0(X)$ by a smaller Banach
space upon which the transfer operator has good spectral estimates.  

\subsection{The Banach spaces}
Let $\{I_j\}$ denote
the arcs in the Bowen--Series coding (not to be confused with intervals from the
definition of~$T$).

\noindent
\begin{minipage}{120mm}
\begin{enumerate}
  \item Choose (small) open neighbourhoods $ I_j \subset U_j \subset \mathbb C$
      such that for any $ I_j \subset T(I_k)$ we have that the closure $
      \overline{ U_j} \subset T(U_k) $;
  \item Let $B \subset C^0( \sqcup_k U_k )$ be the Banach space of bounded analytic functions on the disjoint union 
    $\bigsqcup_k U_k$ with the supremum norm 
    $$
    \|f\|= \sup_{z\in\bigsqcup_k U_k}|f(z)|;
    $$
  \item By construction we see that the Banach space is preserved by all
      $\mathcal L_s$, i.e.,  $\mathcal L_s\colon B \to B$ for all $s\in \mathbb
    C$;
    \setcounter{mylist}{\value{enumi}}
\end{enumerate}
\end{minipage}
\begin{minipage}{3mm} \quad
\end{minipage}
\begin{minipage}{45mm}
\includegraphics{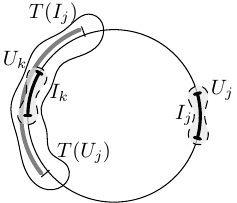}
\captionof{figure}{Arcs~$I_j\subset \uppartial\mathbb D$ of a Markov partition
of~$T$, their images and neighbourhoods. }
\end{minipage}

\begin{enumerate}
  \setcounter{enumi}{\value{mylist}}
  \item On the Banach space $B$ we have that  $\mathcal L_s$ is trace class;
    in other words, there exists only countable many non-zero eigenvalues and
    their sum is finite\footnote{This follows
    directly from the fundamental work of Grothendieck~\cite{G55},~\cite{G56} and
    Ruelle~\cite{R76}. We also refer the reader to~\cite{M91}.}. The finiteness
    of trace allows us to write
    $$
   \det(I- \mathcal L_s) \colon = \exp\left(-\sum\limits_{n=1}^\infty
  \frac1n \trace (\mathcal L_s^n) \right),
    $$
 where the right hand side converges to an analytic function for all $s \in \mathbb C$; and
  \item Furthermore the connection between the derivative of~$T$ and the length
      of closed geodesics allows one to establish 
  $  Z_V(s) =  \det(I- \mathcal L_s )  $; Finally
  \item  We have the classical observation:  
    $Z_V(s) = 0$ if and only if there exists a non-zero function $f \in C^\omega( \sqcup_k U_k)$ such
    that $\mathcal L_s f = f$. 
 \end{enumerate}   
 \noindent We will return to the use of Banach spaces later in the special case of
 infinite area surfaces, where the compactness enters via the limit
 set~$\Lambda$.  

\subsection{Mixing of geodesic flow}
We next consider a closely related  problem in the general area of smooth ergodic theory.
We begin with a definition.

\begin{definition}
    \label{def:corfun}
    Let $Y$ be the unit tangent bundle to a surface~$V$ and let 
    $\nu$ be a (normalised) Liouville measure on~$Y$.
 Given two smooth functions $F,G \colon Y \to \mathbb{R} $ and a
 geodesic  flow $\varphi_t$ on~$Y$, we define a correlation function: 
    $$
    \rho(t) = \int_Y F(\varphi_tx) G(x) d\nu(x) - \int_Y F d \nu \int_Y G d
    \nu, \quad t \ge 0,
    $$
    We say that the flow is strong mixing if $\rho(t)\to0$ as $t\to+\infty$. The rate of convergence in
    this case is referred to as \emph{the rate of mixing}.
\end{definition}
    
In the case of constant negative curvature the Liouville measure corresponds to
the normalised Haar measure.  
    
Mixing is a property in the hierarchy of ergodic properties. It implies
ergodicity and is in turn implied by the Bernoulli property. In the case
of square integrable functions $F, G \in L^2(\nu)$ there is no reason to
expect to be able to say anything about the speed of mixing (i.e., the rate at which
$\rho(t)$ tends to zero). But since we are assuming~$F$ and~$G$ are smooth more
can be said. 
    
Of course the definition can also be reformulated for different invariant
measures, but for our present purposes we only need to consider the measure~$\nu$.  
    
The following result dates back over sixty years and was the forerunner of
method that has proved to be very successful over subsequent years.
\begin{theorem}[Fomin--Gelfand, 1952]
\label{mixing}
The geodesic flow on a manifold of constant negative curvature, 
is strong mixing, in other words $\rho(t) \to 0$ as $t\to + \infty$.
\end{theorem}

The proof of Theorem~\ref{mixing} uses unitary representation of $SL(2,\mathbb R)$ on the space $L^2(Y)$ given by  
$[U_{\mapg} f](x) = f(\mapg^{-1}x ) $.
The representation $U_{\mapg}$ is reducible and therefore 
the problem of establishing strong mixing is reduced to the understanding of the
action of each of the component representations. This brings us to the next question: 
\begin{problem}
Can we estimate the speed of convergence $\rho(t) \to 0$ or improve on these
results in any way? 
\end{problem}
  
 The following results on mixing don't directly rely on properties of the zeta function, but the underlying mechanisms are
similar as will hopefully soon become apparent.  
\begin{theorem} 
  \label{correlation:thm}
 Let~$V$ be a compact surface of constant curvature $\kappa = -1$ and let $Y$ be
 its unit tangent bundle. Given two functions $F, G \in C^\infty(Y)$ 
 there exist $\varepsilon>0$, $C>0$ such that the correlation function of the
 geodesic flow~$\varphi_t$ satisfies 
 $$
      |\rho(t)| \le C e^{-\varepsilon t}, \quad t \ge 0, 
 $$
 where~$C$ depends on~$F$ and~$G$, but~$\varepsilon$ depends only on the geodesic
 flow.
\end{theorem}

The argument is based on representation theory approach. In fact the
correlation functions  correspond to ``decay of matrix coefficients'' in
representation theory~\cite{M87}.
        
We have used the same constant $\varepsilon > 0$ for both results, in the
counting geodesics problem Theorem~\ref{counterror}, and for the speed of mixing
of geodesic flow in Theorem~\ref{correlation:thm}. Assuming this wasn't carelessness,
one might ask the following question: 
        
\begin{problem}
  How are the $\varepsilon$'s related in two problems? 
\end{problem}
The simple answer is:  {\it ``They are the same!''}. To see this, we can
consider the Laplace transform of the correlation function~$\rho$: 
$$
\hat \rho(z) = \int_0^\infty e^{-zt} \rho(t) \mathrm{d} t, \quad z \in \mathbb C.
$$
Thus we are introducing another  complex function to describe the mixing rate.  
The function~$\hat\rho$ is well-defined for $\Re(z) > 0$ since the integral converges on this
half-plane. However, there is much more that one can say.  We  recall the
following properties of~$\hat \rho $.
  
\begin{lemma}
The Laplace transform of the correlation function has the following properties:
\begin{enumerate}
\item The function~$\hat \rho$ extends meromorphically to~$\mathbb C$; 
\item Poles $z_n$ for $\hat \rho$ are related to zeros $s_n$ for
the zeta function $Z_V$ by $z_n = s_n -1$; and 
\item The absence of poles for the Laplace transform $\hat \rho$ for a large half-plane
  $\Re(z)>-\varepsilon$ implies exponential decay of~$\rho$ (here $\varepsilon$
  is the same as in Theorem~\ref{correlation:thm}).
\end{enumerate}  
\end{lemma}
The existence of meromorphic extension of~$\hat \rho$ can be obtained by
establishing a connection with the resolvent of
the transfer operators and using properties of the spectrum of the
transfer operator. The second part follows from the observation that~$Z_V$ has a
similar connection to the spectra of a suitable transfer operator.  
The last part is classical harmonic analysis using the  Paley--Wiener theorem.
    
\begin{definition}
  The poles of the meromorphic extension of~$\hat\rho$ (or sometimes
  zeros of~$Z_V$) are called \emph{resonances}. 
\end{definition}

\section{The Bowen--Series map for infinite area surfaces}

Henceforth, we now want to consider  the case of infinite area surfaces~$V$.
However, we shall only consider the case that~$V$ has no cusps, or equivalently
the corresponding Fuchsian group~$\Gamma$ has no parabolic points.

We want to construct the boundary map associated to the action of the group on the limit set.
Recall that by Lemma~\ref{limdim:lem} in this case the limit set~$\Lambda$ is a Cantor set.  
The same basic approach applies --- except that it is even easier: this time we
have to define an expanding map $T \colon \Lambda \to \Lambda$ on the Cantor set. 
This is best illustrated by an example.  We recall our basic example of  the three funnelled surface.

\medskip

\noindent\begin{minipage}{120mm}
\begin{example}[Three funnelled surface]
    \label{ex:pants2}
Let~$V$ be the three funnelled surface described in Example~\ref{pants:ex}. We can
write~$V = \mathbb D/\Gamma$ where~$\Gamma = \pi_1(V)$ is the associated
Schottky group of isometries. Since~$V$ homotopically equivalent to figure eight, $\Gamma = \langle
\mapg_1, \mapg_2 \rangle$ as a free group on two generators $\mapg_1,\mapg_2 \in \mathop{\mathrm {Isom}}(\mathbb D)$. 
One classical construction is to take two pairs of disjoint geodesics
$\gamma_0,\gamma_1,\gamma_2,\gamma_3$ and isometries identifying geodesics
within each pair and mapping ``exterior'' of $\gamma_j$ into ``interior'' of
$\gamma_{(j+2)\mod\!4}$, as pictured in Figure~\ref{fig:c-set}, so that $\mapg_0 \colon
\gamma_0 \to \gamma_2$, $\mapg_1 \colon \gamma_1 \to \gamma_3$, and $\mapg_2 =
\mapg_0^{-1}$, $\mapg_3 = \mapg_1^{-1}$.
Let~$\Lambda_{j}$ be the part of~$\Lambda$ enclosed between the end points
of~$\gamma_j$. Thus we have a  partition $\Lambda =\sqcup_j \Lambda_j$. 
\end{example}
\end{minipage}%
\begin{minipage}{5mm}\quad \end{minipage} %
\begin{minipage}{44mm}
    \includegraphics{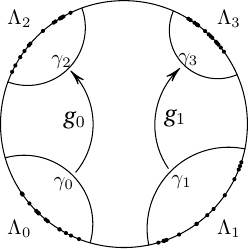}
    \captionof{figure}{Partition of the Cantor set~$\sqcup_j \Lambda_j$.} 
\label{fig:c-set}
\end{minipage}    

\medskip

We can then define a Bowen--Series map $T: \Lambda \to \Lambda$ by $T(x) = \mapg_j(x)$ for $x
\in \Lambda_j$. 
It is easy to see that this is strictly expanding (i.e., $|T'(x)| > 1$ for $x\in
\Lambda$) and Markov, in the sense that the image $T(\Lambda_{j })$ is the union of
the three elements of the partition
$$
\Lambda_{j} = \bigcup_{ k \ne (j+2)\mod\!4 }  \Lambda_k. 
$$

In Lemma~\ref{limdim:lem} we introduced Hausdorff dimension~$\delta$ of the limit
$\Lambda \subset \uppartial\mathbb D$. It turns out we can write~$\delta$ in
terms of the transformation $T \colon \Lambda\to\Lambda$ and its derivative. 
In order to formulate this we need the following definition of a pressure functional. 

\begin{definition}
  \label{pressure:def}  
We define the {\it pressure functional}  $P: C^0(\Lambda)\to \mathbb R$  by 
$$ 
P(f) \colon = \lim_{n\to\infty} \frac1n \log \left(  \sum_{T^nx=x} \exp\left(\sum_{j=0}^{n-1} f(T^j x)\right)\right).
$$   
Equivalently, we can define the pressure using the variational principle:
$$
P(f) \colon = \sup \left\{h_{\mathrm{top}}(\mu) + \int f d \mu \mid \mu
\mbox{ is a } T\mbox{-invariant probability measure}.
    \right\} 
$$
\end{definition}
It is a classical result that in the present context of expanding maps the limit
always exists and finite, we refer the reader to~\cite{PP90} for details. 
We now have the following promised characterisation of $\delta=\mathop{\mathrm
{dim}}_H \Lambda$ in terms of the functions on the boundary. 

\begin{lemma}
  The Hausdorff dimension $\delta= \mathord{\mathrm{dim}}_H \Lambda$  satisfies
  $P(-\delta \log |T^\prime(x)|)=0$. 
 
\begin{figure}[h]
\centering
\begin{tikzpicture}[scale=0.32]
\draw [black, ->] (0,0) -- (16,0);
\node [below] at (16,-0.5) {$t$};
\draw [black, ->] (0,0) -- (0,16);
\node [left] at (13,7) {$P(- \log |T^\prime(x)| \cdot t )$};
\draw (-1,14) .. controls (2,3) and (5,1) .. (12,-1);
\node [below] at (9,0) {$\delta$};
\end{tikzpicture}
\caption{A typical plot of the pressure function.}
\label{pressure:fig}
\end{figure}
\end{lemma}
 
By analogy with the case of compact surfaces, we can introduce  a  family of
complex Ruelle--Perron--Frobenius transfer operators $(s\in\mathbb C) $
initially defined on continuous\footnote{in the induced topology on~$\uppartial\mathbb D$.}
 functions~$C^0(\Lambda)$ on the limit set~$\Lambda$.  

\begin{definition}
We can define the {\it transfer operators} $\mathcal L_s \colon C^0(\Lambda) \to
C^0(\Lambda)$ by 
\begin{equation}
  \label{OpLs:def}
\mathcal L_s f(x) = \sum_{Ty=x} e^{-s \log|T^\prime(y)|} f(y), \quad
s\in \mathbb C.
\end{equation}
\end{definition}

\noindent This gives us another interpretation of the value~$\delta =
\dim_H\Lambda$.

\begin{lemma}
The parameter value $s = \delta$ corresponds to the case when the
operator~$\mathcal L_s$ has the maximal eigenvalue~$1$.
\end{lemma}

As before, to study the zeta function $Z_V$  using the transfer operators we need to
replace $C^0(\Lambda)$ by a smaller Banach space upon which the transfer
operator has good spectral estimates.  

\section{Dynamical applications:  non-compact surfaces}
 
We want to find analogues of the counting and mixing results from
Theorems~\ref{counterror} and~\ref{correlation:thm}, respectively,  in the case
that we have a surface of infinite area. We begin with the counting problem for
closed geodesics.  

\subsection{Counting closed geodesics}
As in the case of compact surfaces, one has in the non-compact case that there
are infinitely many closed geodesics (with one in each free homotopy class). The
obvious question is the following. 
\begin{problem}
What are the analogous results on the asymptotic behaviour of the number of
closed geodesics  in non-compact case? 
\end{problem}
Assume that $V = \mathbb D/\Gamma$ is a surface of infinite area,
where~$\Gamma$ denotes the associated Fuchsian group. We first recall 
that~$\delta$ denotes the Hausdorff dimension of the limit set and let us again
denote by $\# \{ \gamma \colon  \ell(\gamma) \le t \} = \colon N(t)$ the number of
prime closed geodesics of length at most~$t$.  Now we can state the first
result, which gives that~$\delta$ is also the asymptotic growth rate of~$N(t)$.   
\begin{lemma}
The value~$\delta = \dim_H \Lambda $ is the growth rate of the number of
closed geodesics:   
$$
\delta = \lim_{t \to \infty} \frac{1}{t} \log \# \{ \mbox{prime closed geodesic
} \gamma \mbox{ with }  \ell(\gamma) \le t \}.
$$
\end{lemma}
\noindent   
The stronger asymptotic  result in the case of closed geodesics for compact
surfaces (Theorem~\ref{huber}) can also be extended to this case. 
\begin{theorem} \label{infiniteasymp}
We can write 
$$
N(t) \sim \frac{e^{\delta t}}{\delta t} \mbox{ as } t \to +\infty
$$
\end{theorem}
\noindent(Compare with the compact case result in Theorem~\ref{huber} taking into account
that in that case the limit set is the entire circle and $\delta=1$.) 
   
However, whereas the proof in the case of compact surfaces used the Selberg
trace formulae this is not available to us in the case of general infinite area
surfaces. Instead in the present setting the transfer operators can be used
to give an extension to the zeta function~$Z_V$ and this leads to a proof of
Theorem~\ref{infiniteasymp}.   
Similarly to the compact case there is a stronger form of the asymptotic formula with an error term.
\begin{theorem}[Naud 2005]
\label{Naud2005}
There exists $\varepsilon>0$ such that
$$
N(t) = \int_2^{e^{\delta t}} \frac{1}{\log u} \mathrm{d} u \left(1+O(e^{-\varepsilon t})\right) =
\mathop{\mathrm{li}}(e^{\delta t}) + O(e^{(\delta-\varepsilon) t}), 
\mbox{ as } t \to + \infty.
$$
\end{theorem}
        
 Both results Theorem~\ref{Naud2005} and Theorem~\ref{infiniteasymp} require 
 the dynamical approach of transfer  operators. This method also applicable in the case of compact surfaces, but  in that
 case we had the luxury of using the trace formulae, which served us even better
 in as much as they gave both an analytic extension and additional information
 on the location of the zeros. Despite its apparent greater generality, or
 perhaps because of it, the transfer operator approach has the downside that
 there is less control over the location of zeros. 
 The next lemma is the analogue of Lemma~\ref{closed} for compact surfaces.
 
 \begin{lemma} 
 \label{closed2}
 There is a one-to-one correspondence between
 \begin{enumerate}
 \item Closed (prime) geodesics $\gamma$; and 
 \item Periodic (prime) orbits of the Bowen--Series map $ \{x,Tx,\ldots,T^{n-1}
     x \}$ where $T^n x=x $.
 \end{enumerate}
 \end{lemma}
\noindent 
In fact, Lemma~\ref{closed2} has one improvement on Lemma~\ref{closed} in as
much as there is a bijection, without having to concern ourselves with a finite
number of exceptional geodesics. 

As in the  case of compact surfaces we have the following correspondence for the
lengths (compare with Lemma~\ref{closed1}). 
       
\begin{lemma}
To every prime closed geodesics~$\gamma$ which is not a bounding curve
we can associate a periodic (prime) orbit of a Bowen--Series map $\{x,Tx,\ldots,T^{n-1} x \}$,  $T^n x=x$ where $\ell(\gamma)
= \log|(T^n)^\prime(x) |$.
\end{lemma}
We again define the Selberg zeta function
\begin{equation}
  \label{zetaInf:def}
Z_V(s) = \prod_{n=0}^\infty \prod_\gamma \left( 1- e^{-(s+n)\ell(\gamma)} \right)
\end{equation}
where the product is taken over all prime closed geodesics~$\gamma$ and
$\ell(\gamma)$ is the length; as before in~\eqref{zeta:def}. Note
that in the case of infinite area surfaces the closed geodesics are contained 
in the {\it recurrent} part of the flow. The infinite product on the right hand
side of~\eqref{zetaInf:def} converges for $\Re(s)>\delta$ as is easily seen
using Theorem~\ref{infiniteasymp} and defines an analytic function. 

To proceed we again replace $C^0(X)$ by a smaller Banach space upon which the
transfer operator has good spectral estimates   
(e.g. H\"older, $C^1$, $C^k$, $C^\omega$,
the smaller space the better the results).

\subsection{The Banach spaces}
The main case we want to concentrate on is the three funnelled surface (or
``pair of pants'') from Example~\ref{ex:pants2}. This will be a paradigm for a
general case. 

In order to study the  zeta function~$Z_V$ using the transfer operator we need 
to be more restrictive in the  type of  Banach space upon which the transfer operator acts.
We will consider a Banach space of analytic functions.
\begin{enumerate}
  \item Choose four disjoint open neighbourhoods $U_j \supset \Lambda_{j }$ of
      elements of partition of the Cantor set, as pictured in
      Figure~\ref{fig:c-set};
\item Let $B \subset C^\omega( \sqcup_j U_j )$ be the \emph{Banach space} of bounded analytic
functions on the disjoint union of the disks $\sqcup_{j} U_{j }$ with the norm
$$
\|f\| = \sup_{z \in \sqcup_{j} U_{j }} |f(z)|.
$$
Of course, this is equivalent to the direct sum of spaces of analytic functions on
each of the four neighbourhoods.
\item By construction we see that the Banach space is preserved by transfer
  operators given by~\eqref{OpLs:def}, i.e., ${\mathcal L}_s \colon B \to B$ for
  every $s \in \mathbb C$.
  
\item On the Banach space~$B$ we have that~$\mathcal L_s$ is trace class,
   in particular, there exists only a countable set of non-zero eigenvalues which sum is
  finite\footnote{This again follows directly from the works of 
  Grothendieck~\cite{G55},~\cite{G56} and Ruelle~\cite{R76}.}. Thus we can write 
$$
\det(I- \mathcal L_s) \colon = \exp\left(-\sum\limits_{n=1}^\infty
\frac 1n \trace (\mathcal L_s^n) \right)
$$
and the right hand side converges to a function analytic on $ \mathbb C$. 
\item Furthermore using the connection between the derivatives of the
    Bowen--Series map and closed geodesics one can show that $ Z_V(s) =  \det(I- \mathcal L_s) $.
\item Finally, we have the classical observation: $Z_V(s)=0$ if and only if
  there exists a function $f\ne0 \in B$ such that $\mathcal L_s f = f$. 
Thus in order to show that there is a zero free strip it suffices to show that 
there exists $\delta^\prime < \delta$ such that for any $\delta^\prime < \Re(s) <
\delta$ the operator~$\mathcal L_s$ doesn't have~$1$ as an eigenvalue.
 
\end{enumerate}   

As we have observed already,  closed geodesics in the case of non-compact
surfaces lie in the recurrent part of the geodesic flow. We also adopt the
convection that they are oriented geodesics (i.e., a pair of geodesics which are
identical as sets but have different orientations counted as two different geodesics). 

We recall basic properties of the zeta function (cf.~\cite{P91}) for infinite
area surfaces.
\begin{lemma}
  \label{zetaMain:lem}
Let $V$ be an infinite area surface of constant negative Gaussian curvature~$-1$.  
\begin{enumerate}
  \item The infinite product~\eqref{zetaInf:def} converges for $\Re(s)>\delta$
    and therefore $Z_V$ is a well defined analytic function on this half-plane;
\item The zeta function $Z_V$ has an analytic extension to~$\mathbb C$.
\end{enumerate}
\end{lemma}
   
 
 
\subsection{Measures and mixing}
To extend the theory to the case of surfaces of infinite area we first need a
new measure replacing the Haar measure. In this case we are concerned with
measures supported on the recurrent  part of the geodesic flow. The easiest way
to construct such measures is by using the classical construction of measures on
the limit set and then to convert these into measures invariant by the geodesic
flow and supported on the recurrent part of the flow.
 
\begin{definition}
   Let  $\delta_{\mapg 0}$ is the Dirac delta measure supported on the
   point $\mapg 0$.
   We can define a measure on the limit set $\Lambda$ by 
   $$
   \nu_\delta(A) = \lim_{t \searrow \delta} \frac{1}{t - \delta} \cdot 
   \frac{\sum_{\mapg  \in \Gamma} \delta_{\mapg 0}(A) e^{-t d(0,\mapg 0)}}
   {\sum_{\mapg  \in \Gamma} e^{-\delta d(0,\mapg 0)}},
   \qquad A \subset \Lambda.
   $$
   \end{definition}
  
   This is  the standard construction of the Patterson--Sullivan measure on~$\Lambda$
   (cf. a book of Nicholls~\cite[Ch.~3]{Nicholls}, for example).
   Although for each parameter~$t$ we obtain a measure on the countable set
   of points in the orbit of the point~$x_0$, and these measure live on the
   open disk and give zero measure to the boundary circle (and the limit set). 
   As~$t$ decreases more weight is given to points closer to the boundary (in the Euclidean
   metric). In the limit these measures converge in the weak-$*$ topology of
   measures (on the closed unit disk) to a measure actually supported on the
   boundary cicle. 

   We use the Patterson--Sullivan  measure $\nu_\delta$ and the Lebesgue measure
   along the flow direction~$d t$ to define a flow invariant measure 
   $\widehat \nu_\delta$ on the unit tangent bundle to the universal cover,
   which can be identified with $\bigl(\uppartial\mathbb D \times \uppartial\mathbb D
   \setminus \mbox{diagonal} \bigr) \times \mathbb R$ by:  
   $$
   d \widehat \nu_\delta (x,y,t) = \frac{d \nu_\delta \times d \nu_\delta \times d
   t}{|x-y|^{2\delta}}
   $$
   This    corresponds to an invariant measure on $\uppartial\mathbb D \times
   \uppartial \mathbb D$ under the
   diagonal action of $\Gamma$, i.e., $\mapg (x,y) = (\mapg x,\mapg y)$.  Finally, by considering the quotient by $\Gamma$ and
   normalizing we
   arrive at a~$\varphi_t$-invariant measure~$\nu$ on the unit tangent bundle~$Y$. 

  Similarly to the compact case (Definition~\ref{def:corfun}), given two smooth functions $F,G \colon Y \to \mathbb R$ we define the
   correlation function associated to the flow $\varphi_t$ by 
   $$
\rho(t) = \int_Y (F \circ \varphi_t ) G d \nu - \int_Y F d \nu \int_Y G d \nu 
$$ 
In the case of infinite area surfaces the approach of unitary representations
does not apply so naturally.  However, there are dynamical approaches, mostly
rooted in the work of Dolgopyat~\cite{D98}, that gives an analogue to
Theorem~\ref{correlation:thm}.

\begin{theorem}[Naud, 2005] There exists $C, \varepsilon>0$ such that $|\rho(t)|
  \le C e^{-\varepsilon t}$ for $t \ge 0$. 
\end{theorem}
The proof makes use of transfer operators.  These are used to  analyze  the
Laplace transform~$\hat \rho$
by showing that the resolvent of the transfer 
operator~$\mathcal L_s$ given by~\eqref{OpLs:def} has no poles for $\Re(s) \ge -\varepsilon$. 

\section{Location of resonances for infinite area surfaces}
\label{s:locres}
Finally, we want to discuss the behaviour of the distribution of the zeros for
the analytic extension of the Selberg zeta function in the case of infinite area surfaces.  

\subsection{Rigorous results: the case of pair of pants}

To demonstrate the method, we will consider the case of a pair of pants. This is
the simplest example to study and in this case the results are clearer to
interpret.  However, the analysis can be adapted to other examples and
some of these are discussed in the last section.

\subsubsection{Distribution of zeros}
Since the Hausdorff dimension corresponds to the parameter value when the
operator $\mathcal L_s$ has maximal eigenvalue~$1$, and the zeta function is
$\det(I - \mathcal L_s)$, we may deduce that $Z_V(\delta) = 0$. 
We also know~\cite{Naud1} that there are no zeros in a (possibly very) small strip
$\delta^\prime < \Re(s) < \delta$.    A natural question so ask is what happens when we first venture into the region
$\Re(s) < \delta^\prime$?  
\begin{problem}
  What about the other zeros not so near $\Re(s)=\delta$?
\end{problem}

Since we expect that there will be zeros in a larger region, it is therefore
more appropriate to ask about the density of zeros. Naud showed the following:
\begin{theorem}[Naud, 2014]
  Fix $\sigma \in \left(\frac{\delta}{2},\delta\right)$. Then there exists
  $0<\eta < \delta$: 
  $$
   \# \{ s_n = \sigma_n + it_n \mid Z_V(s_n)=0  \colon \sigma \le \sigma_n \le \delta,\,|t_n|\le t
   \} = O\left(t^{1+\eta}\right).
  $$
\end{theorem}
The proof is particularly elegant and uses both the Bowen--Series map and  properties of the  pressure.

\subsubsection{Patterns of zeros}
Our final result concerns the strange patterns of zeros that one sees empirically.
Consider again the example of a pair of pants and associate the Selberg zeta
function (cf.~\eqref{zetaInf:def})
    $$
    Z_V(s) = \prod_{n=0}^\infty \prod_\gamma
    \left(1-e^{-(s+n)\ell(\gamma)}\right), \qquad s \in \mathbb C,
    $$ 
where~$\gamma$ is a (prime) closed geodesic or an orbit of the flow of length~$\ell(\gamma)$. 
We have already mentioned (cf. Lemma~\ref{zetaMain:lem}) that the infinite product
converges for $\Re(s) > \delta = \dim_H \Lambda$, the Hausdorff dimension of the
limit set of the associated Bowen--Series map and there is an analytic extension
to~$\mathbb C$, that can be established using trace of the associated transfer operator. 

The main object of study in this section is the zero set of the function~$Z_V$: 
\begin{equation}
  \label{eq:zeroset}
 \mathcal S_V\colon = \{s \in \bbC \mid Z_V(s)=0\}. 
\end{equation}
\begin{problem}
  Where are the zeros of $Z_V(s)$, in strip $0 < \Re(s) < \delta$ located?
\end{problem}

\begin{figure}[h!]
\centerline{\includegraphics[scale=0.80,angle=0 ]{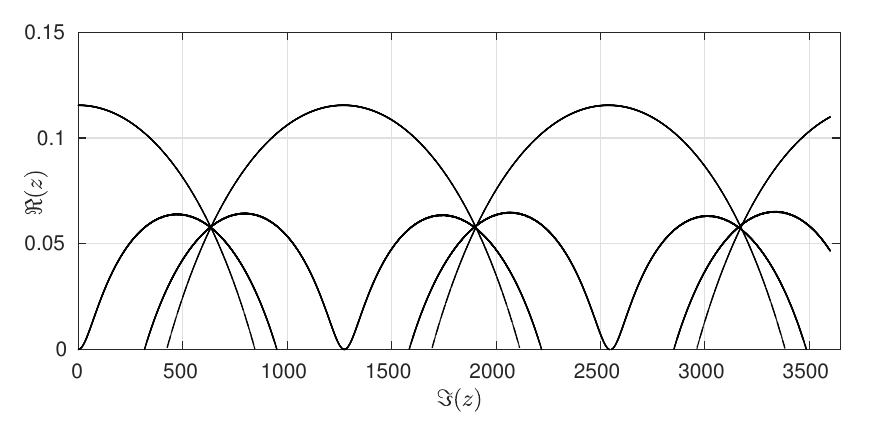} }
\caption[A pattern of zeros]{
 The zeros of the associated zeta function 
$Z_V$ in the critical strip. 
The individual zeros are so close in the plot that it creates the illusion that
they form well defined smooth curves. } 
\label{fig:draft2}
\end{figure}

In pioneering experimental work, D.~Borthwick has studied the location of the
zeros for the zeta function in  specific examples  of infinite area surfaces~\cite{B14}. 
The plot in Figure~\ref{fig:draft2} is fairly typical for zeros in the critical strip 
for a symmetric three funnelled surface where each of the three simple closed geodesics 
corresponding to a funnel has the same, sufficiently large, length. 
A modern laptop allows one to study symmetric three funneled surface with the length of
the three defining closed geodesics at least~$8$ without much difficulty.

Based on numerical evidence, we make the following simple empirical observations 

\bigskip

\noindent{\bf Informal Qualitative Observations.}
Let $V = V(b)$ be the  three funnelled surface defined by three simple closed
geodesics of equal length\footnote{This normalization makes formulae in subsequent
calculations shorter.}~$2b$ then for a sufficiently large $b$: 
\begin{description}
\item[O1]: The set of zeros $\mathcal S_{V}$ {\it appears} to be an almost periodic set 
with translations 
$$
\tau = \{i( \pi k e^b + \varepsilon_k) \mid k \in \mathbb N \}
$$ 
 where $\varepsilon_k = \varepsilon_k(b) = o(e^{-b/2})$ as $b \to +\infty$ for each~$k$.
\item[O2]:
The set of zeros $\mathcal S_{V}$ {\it appears} to lie on a few distinct curves, which {\it seem} to have
a common intersection point at $\frac{\delta}{2} + i \frac{\pi}{2}e^b$, as $b \to +\infty$.
\end{description}

One way to understanding these observations is to find an approximation of the zeta 
function by a simpler expressions whose zero set can be described. 
The first obstacle is that for large~$b$ the value of~$\delta=\mathord{\mathrm{dim}}_H \Lambda$ will be small.
We begin by recalling  the following result of  McMullen~\cite{M98}:

\begin{lemma} 
  The largest real zero $\delta\sim \frac{\ln 2}{b}$ as $b \to +\infty$, i.e.,
  $\lim\limits_{b \to +\infty} b \delta  = \ln 2$.
\end{lemma}
\noindent
It gives us the width of the critical strip: we see immediately that in the limit $b\to+\infty$
the zero set ``converges'' to the imaginary axis. We therefore apply an affine
rescaling that allows us to see the pattern in the zero set for large values of~$b$. A natural choice for rescaling 
factors is the approximate period of the pattern in the imaginary direction and 
the approximate reciprocal of the width of the critical strip in the real direction. 

\begin{notation}
 We will be using the following notation.
  \begin{enumerate}
\item[(a)]
A compact part of the  critical strip of height $h>0$ which we denote by 
$$
{\mathcal R}(h) = \{ s \in \mathbb C \mid 0 \leq  \Re(s)  \leq
\delta \mbox{ and }   |\Im(s)| \leq h\};
$$
and a compact part of the normalized critical strip of height $\hat h>0$ which we denote by 
$$
\widehat{\mathcal R}(\hat h) = \{ s \in \mathbb C \mid 0 \leq  \Re(s)  \leq
\ln 2 \mbox{ and }   |\Im(s)| \leq \hat h\}.
$$
\item[(b)] We denote the rescaled set of zeros by $$ \widehat{\mathcal S}_{V} \colon = 
\left\{\sigma b+i e^{-b} t \,\Bigl| \, \sigma+it \in {\mathcal S}_V \right\}.
 $$
 where evidently, $0< \Re(\widehat{\mathcal S}_{V}) \le \ln 2$.
\end{enumerate}
\begin{figure}[h]
\centering
\begin{tikzpicture}[scale=0.32]
\draw (0,-5) -- (1,-5) -- (1,5) -- (0, 5) -- (0,-5);
\draw (10,-5) -- (12.5,-5) -- (12.5,5) -- (10, 5) -- (10,-5);
\draw [<->, dashed] (0,-5.5) -- (1,-5.5);
\node [below] at (0.5,-5.5) {$\delta$};
\draw [<->, dashed] (-0.5,-5) -- (-0.5,5);
\node [left] at (-0.5,0.0) {$h$};
\draw [<->,dashed] (10,-5.5) -- (12.5,-5.5);
\node [below] at (11.25,-6) {$\ln 2$};
\draw [<->, dashed] (9.5,-5) -- (9.5,5);
\node [right] at (8,0.0) {$\hat h$};
\draw[->] (0,-5) -- (4,-5); 
\draw[->] (0,-5) -- (0,7); 
\draw[->] (10,-5) -- (14,-5); 
\draw[->] (10,-5) -- (10,7); 
\node [right] at (1,2.5) {$\mathcal R(h)$};
\node [right] at (12.5,2.5) {$\widehat{\mathcal R}(\hat h)$};
\node [left] at (0,-5.5) {$0$};
\node [left] at (10,-5.5) {$0$};
\end{tikzpicture}
\hskip 2cm
\begin{tikzpicture}[scale=0.4]
\draw (0,-5) -- (1,-5) -- (1,5) -- (0, 5) -- (0,-5);
\draw (4,-2) -- (6.5,-2) -- (6.5,2) -- (4, 2) -- (4,-2);
\draw [<->, dashed] (0,-5.5) -- (1,-5.5);
\node [below] at (0.5,-5.5) {$\delta$};
\draw [<->, dashed] (-0.5,-5) -- (-0.5,5);
\node [left] at (-0.5,0.0) {$\pi e^b$};
\draw [<->,dashed] (4,-2.5) -- (6.5,-2.5);
\node [below] at (5.25,-3) {$\ln 2$};
\draw [<->,dashed] (7,-2) -- (7,2);
\node [right] at (7,0.0) {$\pi$};
\draw [->] (1.5,0) -- (3.5,0);
\node [above] at (2.5,0.0) {$A_b$};
\end{tikzpicture}
\caption{Left:~The strips ${\mathcal R}(h)$ and 
$\widehat {\mathcal R}(\hat h)$; Right: By renormalizing the strip ${\mathcal R}(\pi e^b)$ to 
$\widehat {\mathcal R}(\pi)$ we can compare the zeros of zeta functions for different $b$, as $b$ tends to infinity.
}
\end{figure}

\noindent{\rm We now introduce a family of four curves approximating $\widehat{\mathcal
S}_{V}$ as $b \to +\infty$.  }
\begin{enumerate}
\item[(c)] Let $\mathcal C=\cup_{j=1}^4 \mathcal C_j$, where the $\mathcal C_i$ are explicit curves given by 
$$
\begin{aligned}
\mathcal C_1 &= \left\{ 
\frac12\ln|2-2\cos(t)| + it \, \Bigl| \, \frac\pi3 \le t \le \pi \right\}; \cr 
\mathcal C_2 &= 
\left\{ \frac12 \ln|2+2\cos(t)| + it \, \Bigl| \, 0 \le t \le \frac{2\pi}3 \right\}; \cr
\mathcal C_3 &= \left\{ \frac12 \ln \Bigl|
1 - \frac12 e^{2it} - \frac12 e^{it} \sqrt{4 - 3 e^{2i t}} \Bigr| + it \, \Bigl| \, 0
\le t \le \frac{3\pi}4  \right\}; \cr
\mathcal C_4 &= \left\{  \frac12 \ln \Bigl|
1 - \frac12 e^{2it} + \frac12 e^{it} \sqrt{4 - 3 e^{2i t}} \Bigr| + it \, \Bigl| \,
\frac\pi4 \le t \le \pi \right\}.
\end{aligned}
$$

The key thing to 
note is that the curve $\mathcal C$ in Figure~\ref{fig:c-curves} look very  similar to empirical plots in
Figure~\ref{fig:uni}. 

\end{enumerate}
\end{notation}

\begin{figure}[h]
\centering
\includegraphics{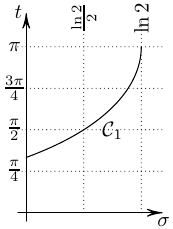} \hskip 0.25cm
\includegraphics{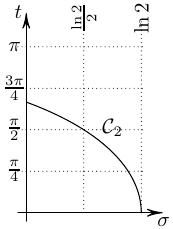}  \hskip 0.25cm
\includegraphics{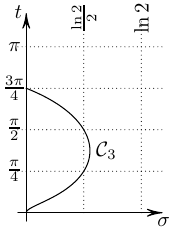}  \hskip 0.25cm
\includegraphics{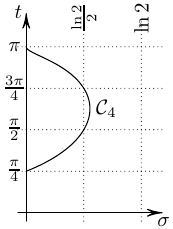}  \hskip 0.25cm
\includegraphics{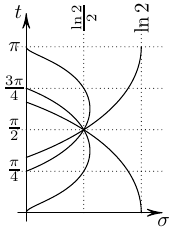} 
\caption{Plots of the curves $\mathcal C_1$,
 $\mathcal C_2$, $\mathcal C_3$, $\mathcal C_4$ and their union $\mathcal C$.}
 \label{fig:c-curves}
\end{figure}

We can now formally state the approximation result, which provides an
explanation for the Observations.
\begin{theorem}
\label{curves}
The sets $\widehat{\mathcal S}_{V}$ and $\mathcal C$ are close in
the  Hausdorff metric~$d_H$ on a large part of the strip $0<\Re(s)<\log 2$.  More
precisely, there exists $\varkappa  > 0$  such that 
$$
   d_H( \widehat {\mathcal S}_{V} \cap \widehat{\mathcal R}(e^{\varkappa
   b}) ,\mathcal C \cap \widehat{\mathcal R}(e^{\varkappa b}) ) \to 0 \mbox{
 as } b \to +\infty.
$$ 
\end{theorem}
  The theorem implies that every rescaled zero $s\in\widehat {\mathcal S}_{V(b)}
  \cap \widehat{\mathcal R}(e^{\varkappa b})$ belongs to a neighbourhood of
  $\mathcal C$ which is shrinking as $b \to \infty$. On the other hand, the rescaled 
  zeros are so close, that the union of their shrinking neighbourhoods contains $\mathcal
  C$. 
  The most significant feature of this result is that the height $e^{\varkappa b}$ of
  the rescaled strip $\widehat{\mathcal R}(e^{\varkappa b})$ is larger than the
  period of the curves $\mathcal C$,
  and it corresponds to a part of the original strip of the height
  $e^{(1+\varkappa)b}$. 

  Because of the natural symmetries of~$V(b)$, it is convenient  to choose a
presentation of the associated Fuchsian group in terms of three  reflections (as
in~\cite{M98}, for example). More precisely, we can fix a value
$0<\alpha\le\frac{\pi}{3}$ and consider the Fuchsian group $\Gamma =
\Gamma_\alpha  := \langle R_1, R_2, R_3\rangle$ generated by reflections $R_1,
R_2, R_3$ in  three disjoint equidistant geodesics $\beta_1,\beta_2,\beta_3$.
In the previous description of the Bowen--Series coding we chose our generators
to be orientation preserving. However, to exploit the symmetry of $V$ these
orientation reversing generators are more convenient, although the same general
theory applies as before.

\begin{figure}[h]\label{fundamental1}
 \centering
  \includegraphics[scale=1]{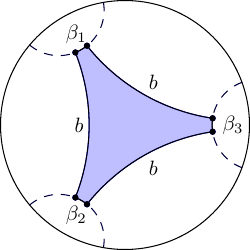}
\caption[Fundamental domain]{Three circles of reflection in Poincar\'e disk with pairwise
distance~$b$.}
\end{figure}
 Although the three individual generators are 
 orientation reversing the resulting quotient  surface $\mathbb D/\Gamma$ is an oriented
  infinite area surface.
  We now explain the pattern of the zeros for $Z_{V}$ (for large~$b$) in terms of those of a simpler approximating function.

\subsubsection{Approximating Selberg zeta function}

Finally, we come to a peculiar issue.
\begin{problem}
  Where do these curves come from? 
\end{problem}
Consider the $6\times 6$ complex matrix function 
$$
B(z) = \left(
\begin{matrix}
  1 & z & 0 & 0 & z^2 & z \\
  z & 1 & z^2 & z & 0 & 0 \\
  0 & 0 & 1 & z & z & z^2 \\
  z^2 & z & z & 1 & 0 & 0 \\
  0 & 0 & z & z^2 & 1 & z \\
  z & z^2 & 0 & 0 & z & 1 \\
\end{matrix}
\right)
$$
We can use this to give an approximation to $Z_{V}$ and its zeros.
\begin{theorem}[Approximation Theorem]\label{approx}
  Using the notation introduced above, the real analytic function $Z_{V}\bigl(\frac{\sigma}{b} +
ite^b\bigr)$ converges uniformly to $\det(I-\exp(-2\sigma-2itbe^b)B(e^{it}))$, i.e., 
$$
\left|Z_{V}\left(\frac{\sigma}b + i
t e^{b}\right) -\det\left(I-\exp(-2\sigma-2itbe^b)B(e^{it})\right)  \right| 
\to 0 
\mbox{ as }  b \to +\infty.  
$$
\end{theorem}

\begin{figure}[h]
\begin{tabular}{cc}
\includegraphics[scale=0.5,angle=0 ]{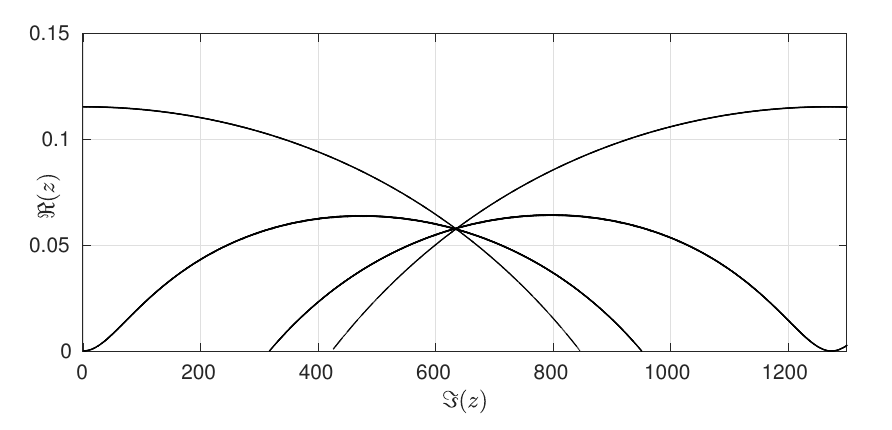} &
\includegraphics[scale=0.5,angle=0 ]{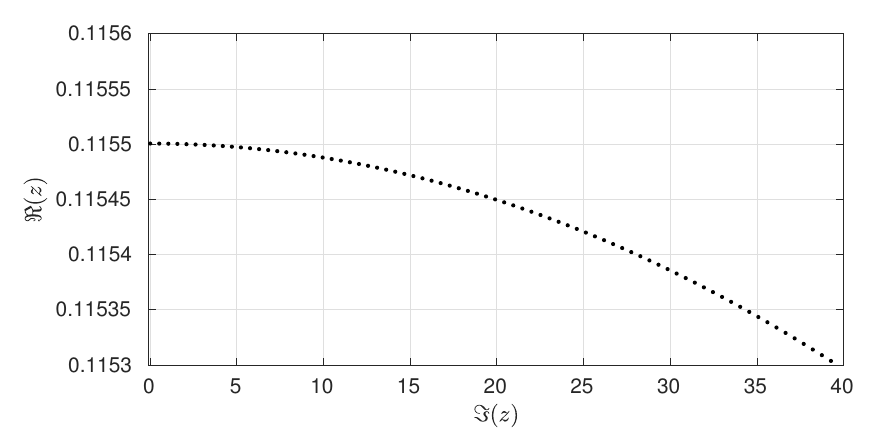} \\
(a) & (b) \\
\end{tabular}
\caption{(a) The zeros of the determinant $\det\left(I -
e^{-2b\sigma} B(\exp(ie^{-b} t))\right)$ approximating $Z_{V(b)}$;  (b) A zoomed
version in a neighbourhood of $\delta$ showing indvidual zeros. The distance
between imaginary parts of consecutive zeros is approximately $\frac{\pi}{b}$.}
\label{fig:uni}
\end{figure}

\paragraph{Key observation:}
The importance of the matrix $B$ is that for $s= \sigma + i t$ to lie on the
curves $\mathcal C_1, \mathcal C_2, \mathcal C_3, \mathcal C_4$
corresponds to eigenvalues $\lambda_t$ for $B(e^{-it})$ to satisfy
$e^{-\sigma}|\lambda_t| =1$. 

\bigskip

We will briefly explain how Theorem~\ref{approx} implies Theorem~\ref{curves}.
We shall show that
for all $\varepsilon>0$ and $T>0$ there exists
$b_0 > 0$ such that for any $b > b_0$ the zeros of the function
$Z\left(\frac{\sigma}{b} + it e^{b}\right)$ with $0 \leq \sigma \leq \ln 2$ and $|t|
\leq e^{(2-\varkappa)b}$ belong to a neighbourhood $\cup_k U(\mathcal C_k,
\varepsilon)$ of the union of the curves $\cup_k \mathcal C_k$.

Given $\varepsilon>0$ and a point $z_0=\sigma_0+it_0$ outside of $\varepsilon$-neighbourhood
  of $\cup_{j=1}^4 \mathcal C_j$ 
  a straightforward computation gives that the determinant
  $$
  \left|\det(I-\exp(-2\sigma_0-it_0be^b)B(\exp(it_0)))\right|>\exp(-6\varepsilon)(\exp\varepsilon-1)^6>0
  $$ 
  is bounded away from zero and the
  bound is independent of $b$. Thus  we see that outside of the neighbourhood $\cup_{j=1}^4 U(\mathcal
  C_j,\varepsilon)$ the determinant has modulus uniformly bounded away from~$0$.  
  But using Theorem \ref{approx} the zeta function
  $Z_{V}\left(\frac{\sigma}b + it e^b\right)$
  can be approximated arbitrarily closely by the determinant.  
  Therefore for $b$ sufficiently large all zeros of the function
  $Z_{V}\left(\frac{\sigma}{b} + it e^b\right)$ belong to the $\varepsilon$-neighbourhood of
  $\cup_{j=1}^4 \mathcal C_j$.


\subsection{The case of one-holed torus}
In this section we want to apply ideas developed in~\cite{PV18} and described
above to study patterns of zeros of the Selberg zeta function associated to a symmetric one-holed torus.
Contrary to the case of a pair of pants, they still remain a mystery, and our
goal here is to formulate a few conjectures on the distributions of zeros,
and to give more insight into our approach. We try to make this section as
self-contained as possible, though we reuse the notation.

The Selberg zeta function corresponding to the holed torus is the same 
as the one associated to a genus one hyperbolic surface with a single funnel, the case studied by
D.~Borthwick and T.~Weich~\cite{B14},~\cite{BW16}.
We start with numerical experiments. 

\subsubsection{Plots of zeros and conjectures}
\begin{figure}[h]
  \begin{center}
  \includegraphics{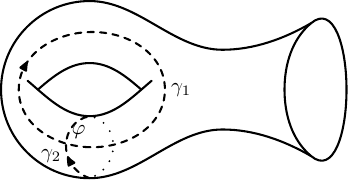}
\end{center}
\caption{One-funneled torus with Fenchel--Nielsen coordinates: a pair of
shortest geodesic $\gamma_1$, $\gamma_2$ and angle inbetween~$\varphi$. Artistic
impression. There is no embedding into $\mathbb R^3$.}
\label{fig:torus}
\end{figure}

It is well-known that, as a Riemann surface, a one-holed torus is uniquely defined by the length of two
shortest geodesics and the angle in between as shown in Figure~\ref{fig:torus}.
We say that a one-holed torus is symmetric if the two geodesics have the same
length and are orthogonal to each other. Therefore a symmetric one-holed torus has
only one parameter and we will denote it by~$\mathring{ \mathbb T}(a)$, where~$a$ 
is the length of the shortest closed geodesic.
We would like to consider three different symmetric tori\footnote{This choice
will be clear later, here we just note that the first choice 
$a=10$ is guided by previous research~\cite{B14} and the integer part $\left[\frac{10}{\log 2}\right] =
14$.} $\mathring{\mathbb T}(10)$, $\mathring{\mathbb T}(14\log2 +0.05)$, and
$\mathring{\mathbb T}(14 \log 2)$.

It is known~\cite{M98} that the width of the critical strip is proportional
to~$a^{-1}$. 
Physical dimensions of paper and screen impose limitations on the figures. To
make them more realistic, we plot rescaled zeros and chose aspect ratio of the
image to correspond the scale on the axis. In addition, the algorithm we are using
allows to compute the zeros in a small part of the critical strip near the real
axis $0<\Re z<\delta$, $0< \Im z <e^{\frac{3a}{2}}$ only. We refer to this
subset of the zero set as ``small zeros'', and these are \emph{the only zeros we
consider here}, unless stated otherwise. 

\begin{figure}
  \caption{Characteristic plots of rescaled zeros of symmetric torus. Larger
  balls mean that the rescaled zeros are further away. In all three subfigures
  $a \cdot e^{-\frac{a}{2}} \approx 0.07$ and we show the strip $0 < \Im \hat z <
  0.2\approx 3\cdot 0.07$.}
  \label{fig:torus-plots}
  \includegraphics{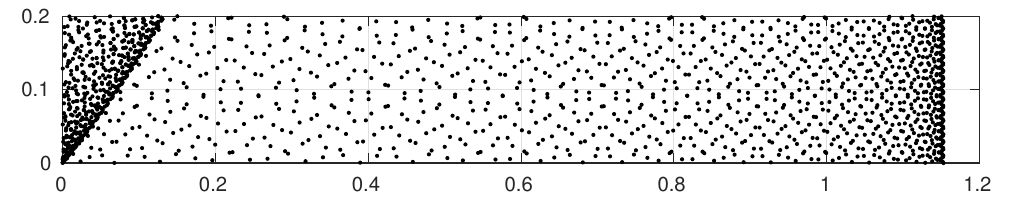}
  \centerline{(a) \ 
  The plot shows $1777$ rescaled zeros $\hat z\in\widehat S(10)$, $0<\Im \hat z <0.2$.}\\
  \includegraphics{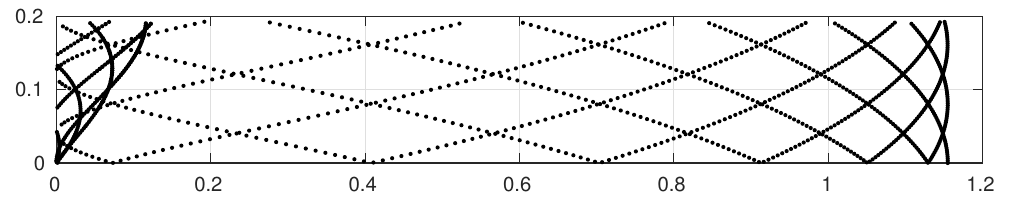}
  \centerline{(b) \  
  The plot shows $1382$ rescaled zeros $\hat z \in \widehat S(14\log2+0.05)$,
  $0<\Im \hat z < 0.2$.} \\
  \includegraphics{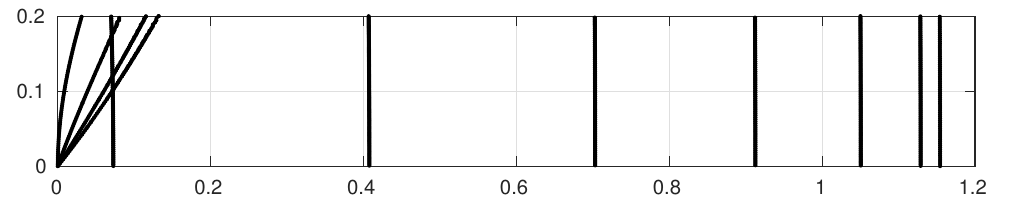}
  \centerline{(c) \ 
  The plot shows $1359$ rescaled zeros $\hat z \in \widehat S(14\log 2)$, $0<\Im
  \hat z < 0.2$.}\\ 
\end{figure}

Three sample plots are shown in
Figure~\ref{fig:torus-plots}. The plot in Figure~\ref{fig:torus-plots}(b) 
corresponds to a small region of the plot in~\cite{BDW17}, p.~30, Fig.~7
(bottom) near the imaginary axis, and the plot in 
Figure~\ref{fig:torus-plots}(c) corresponds to the plot in~\cite{B14}, p.~7, Fig.~11 (left).

We will be using the following notation.  Let 
\begin{equation}
    \widehat{S}(a) = \widehat{S}_{\mathring{\mathbb T}(a)} \eqdef \left\{ \hat z = a \left( \Re z + e^{-2a} \Im z \right) 
\mid Z_{\mathring{\mathbb T}(a)}(z) = 0, \, 0< \Im z<e^{\frac{3a}{2}} \right\}
\label{eq:zhatDef}
\end{equation}
be a set of small rescaled zeros of the zeta function $Z_{\mathring{\mathbb
T}(a)}$ associated to the torus~$\mathring{\mathbb T}(a)$. 

Based on numerical experiments we state several conjectures on
properties of small rescaled zeros, which we will try to explain heuristically
in Section~\S\ref{s:args}, leaving rigorous arguments for another occasion.

\begin{figure}
 \caption{Characteristic plots of rescaled zeros with a cover by disks of the
  radius equal to the distance to the nearest zero. Bigger disks mean that the
  zeros are further apart.}
 
  \includegraphics{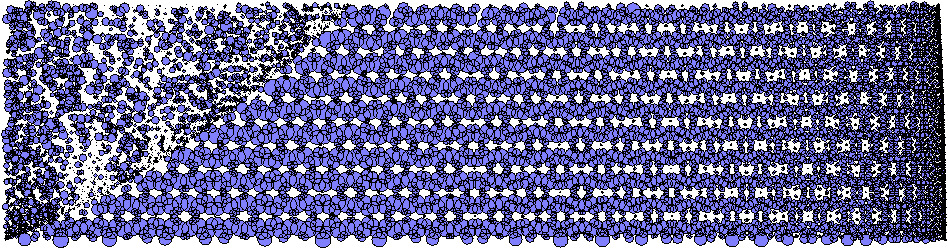}
  \centerline{(a) \ The cover $Cr(\widehat S(10))$ of the zero set from Figure~\ref{fig:torus-plots}(a).}

  \includegraphics{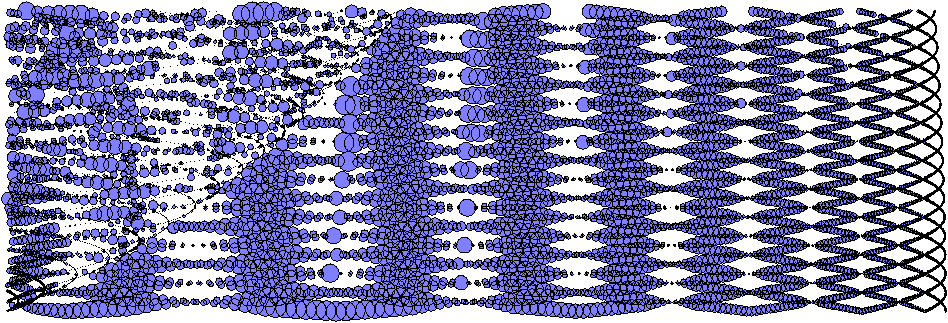}
  \centerline{(b) \ The cover $Cr(\widehat S(14\log 2+0.05))$ of the zero set from Figure~\ref{fig:torus-plots}(b).}
  
  \includegraphics{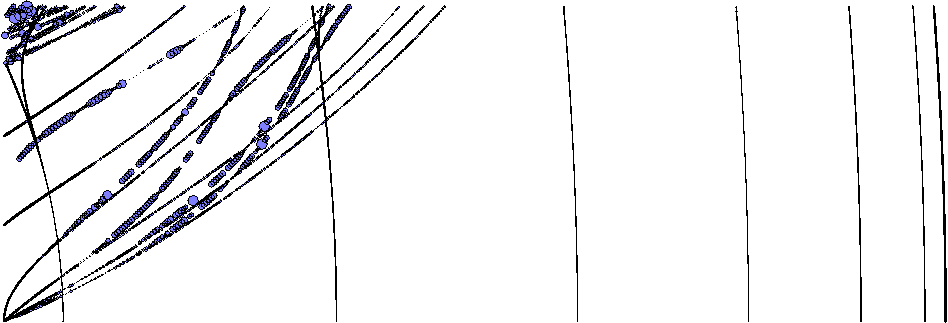}
  \centerline{(c) \ The cover $Cr(\widehat S(14\log 2))$ of the zero set from Figure~\ref{fig:torus-plots}(c).}

  \label{fig:cover}
\end{figure}

To characterise the density of the set of small zeros we
define a cover by open balls of a given set~$A$ by  
\begin{equation}
  Cr(A) = \bigcup_{z \in A} B( z,r(z)), \mbox{ where
  } r(z) =
  \min_{ z'\ne z, z' \in A} | z- z'|.
\label{eq:CrDef}
\end{equation} 
The plots illustrating the cover in the cases which we consider are shown in
Figure~\ref{fig:cover}.

In order to describe differences between the sets of small rescaled zeros corresponding to different tori,
we will study the dependence of the following characteristics on the length
parameter~$a$.
\begin{enumerate}
\item The Hausdorff distance between the convex hull of the set of small
  rescaled zeros and the set itself $\mathord{D}_H (a) \eqdef
  \mathord{\mathrm{dist}}_H(\mathrm{Conv}(\widehat S(a)),\widehat S(a))$;
\item The area of the cover $\mathord{M} (a) \eqdef
  \mathord{\mathrm{Area}}(Cr(\widehat S(a)))$; 
\item Expectation $\mathord{E}(a) \eqdef \mathop{\mathrm E}(\inf|z-\widehat S(a)|)$ and variance 
  $\mathord V(a) = \mathop{\mathrm{Var}}(\inf|z-\widehat S(a)|)$ of the distance from a randomly
  chosen point $z$ in the critical strip to $\widehat S(a)$. 
\end{enumerate}
A selection of three representative empirical results is shown in
Table~\ref{tab:numsum}. On the basis of these (over~$100$ many) 
results we propose the following conjecture. 

\begin{conjecture}[distribution of rescaled zeros near the real axis]
The Hausdorff distance $D_H$, the area function $M$, the expectation $E$, and the
variance $V$ are continuous and monotone functions with respect to the fractional
part $\left\{ \frac{a}{\log 2}\right\}$. In particular, $D_H$ and $M$ are
increasing while $E$ and $V$ are decreasing. Therefore, $a \in \mathbb N \log 2$
are local maxima for $E$ and $V$, and local minima for $D_H$ and $M$.
\end{conjecture}
In the case of $\mathring{\mathbb T}(k \log 2)$ for some $k \in \mathbb N$
 small zeros can be described more precisely.
\begin{conjecture}[the case of rationally-dependent short geodesics]
    Let $a \in \mathbb N \log 2$. 
    Small zeros with $0<\Im z \le e^{\frac{3a}{2}}$ lie on a small number
     of well defined lines. Among them, there are $\frac{a}{\log 2}$ lines
     nearly parallel to the imaginary axis if $\frac{a}{\log 2}$ is odd, 
     and $\frac{a}{2 \log 2}$ lines nearly parallel to the imaginary axis if 
     $\frac{a}{\log 2}$ is even.
\end{conjecture}
The following conjecture addresses properties of the non-rescaled zero set of 
$Z_{\mathring{\mathbb T}(a)}$ independently of number-theoretic properties
of~$a$.
\begin{conjecture}[distribution of zeros on a large scale]
  The zero set of the Selberg zeta function for $\mathring{ \mathbb T}(a)$ has the
following properties
\begin{enumerate}
  \item The asymptotic of the number of zeros in the critical strip with $\Im z<t$ is
$$
\#\{z \in \mathbb C \mid 0<\Re z<\delta, \, 0< \Im z <t, \,
Z_{\mathring{\mathbb T}(a)}(z)=0 \} = \frac{2a}{\pi}t+O(1);
$$
  \item Any rectangle in the critical strip of the height $\frac{\pi}{a}$
    contains at least one zero;
  \item The real parts of zeros are dense in $\left(0,\frac{1}{2}\delta\right)$.
\end{enumerate}
\end{conjecture}

\begin{table}
\caption{Hausdorff distance, expectation, variance, and area of the cover.}
\label{tab:numsum}
\begin{center}
\begin{tabular}{|l|c|lclclcl|}
  \hline
  $a$ & $\bigl\{\frac{a}{\log 2}\bigr\}$ &  $D_H(a)$ & \ & $M(a)$ & \ & $E(a)$ & \ & $V(a)$  \\
  \hline
 $14 \log 2$ & $0$ &  $0.1671\ldots$  & \ & $0.130210\ldots$  & & $0.040\ldots$ & & $0.001\ldots$  \\
 $14 \log 2 + 0.05$ & $0.072\ldots$ & $0.0356\ldots$ & \ & $0.901338\ldots$ & & $0.010\ldots$ & &
  $5.1\ldots\cdot 10^{-5}$ \\
 $10$ & $0.427\ldots$ & $0.0041\ldots$ & \ & $1.060217\ldots$ & & $0.006\ldots$ & &
  $9.6\ldots\cdot10^{-6}$ \\
  \hline
\end{tabular}
\end{center}
\end{table}

\begin{table}
\caption{Dispersion and expectation of the distance to the zero set measured
using different choices of sample points.}
\begin{center}
\begin{tabular}{||p{70mm}|c|c||} 
  \hline
  Sample point distribution & Expectation & Dispersion \\
  \hline
  \multicolumn{3}{||c||}{ Symmetric torus $a=14 \log2$ } \\
  \hline
  regular rectangular $100\times50000$ points & $0.040454\ldots$ & $0.001345\ldots$
  \\
  \hline
  random $50\times1000$ points & $0.040417\ldots$ & $0.001347\ldots$ \\
  \hline
  random $100\times50000$ points & $0.040361\ldots$ & $0.001347\ldots$ \\
  \hline
  \multicolumn{3}{||c||}{ Symmetric torus $a=14 \log2+0.05$ } \\
  \hline
  regular rectangular $100\times50000$ points & $0.010254\ldots$ &
  $5.135236\ldots\cdot 10^{-5}$ \\
  \hline
  random $50\times1000$ points & $0.010212\ldots$ & $5.163147\ldots\cdot10^{-5}$ \\
  \hline
  random $100\times50000$ points &  $0.010218\ldots$ & $5.148808\ldots\cdot10^{-5}$
  \\
  \hline
  \multicolumn{3}{||c||}{ Symmetric torus $a=10 = 14\log2 + 0.29\ldots$ } \\
  \hline
  regular rectangular $100\times50000$ points & $0.005994\ldots$ &
  $9.705229\ldots\cdot 10^{-6}$  \\
  \hline
  random $50\times1000$ points & $0.006002\ldots$ & $9.641849\ldots\cdot 10^{-6}$ \\
  \hline
  random $100\times50000$ points & $0.005972\ldots$ &
  $9.657625\ldots\cdot10^{-6}$
  \\
  \hline
\end{tabular}
\end{center}
\label{table-exp-dist}
\end{table}

\subsubsection{Geometry of a one-holed torus}
\label{s:torus-geometry}
\begin{figure}
  \begin{tabular}{ccc}
\includegraphics{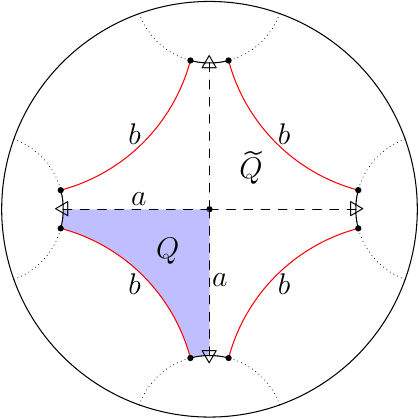} & \qquad \qquad &
\includegraphics{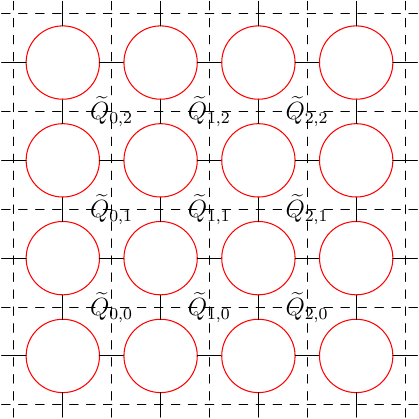} \\ (a) & \ & (b)
  \end{tabular}
  \caption{(a) A right-angled hyperbolic pentagon $Q$ and its three copies
  forming a fundamental domain $\widetilde Q$ in the Poincar\'e disk. (b) A
  tessellation of $\mathbb
  Z \times \mathbb Z$-holed real plane by copies of the fundamental domain
  $\widetilde Q_{i,j}$.}
  \label{pent:fig}
\end{figure}

A very good exposition can be found in Buser and Semmler~\cite{BS88}. Here we
summarise the results we need making necessary adaptations to the case we consider. 
A one holed torus is a genus one Riemann surface whose boundary consists of a
single simple closed geodesic. 
In order to estimate the length of the closed geodesics 
we use Fenchel---Nielsen coordinates and a universal cover by a holed plane.
It turns out that in the case of symmetric one-holed torus the universal cover
has one parameter. Namely, we may consider a right-angled hyperbolic pentagon with the sides
$\{a,*,b,*,a\}$ as shown by the shaded area in Figure~\ref{pent:fig}(a). It is
known cf.~\cite{B83}, \S7.18 that $a$ and
$b$ satisfy the identity
$$
\sinh^2 a = \cosh b.
$$
Fixing the length of the boundary geodesic $b$ we compute
\begin{equation}
a = \frac12 \ln \left( e^{b}  + e^{-b} + 1 + \sqrt{e^{2b}+e^{-2b} + 2(e^b + e^{-b})
+ 1} \right)=\frac12\left(b + \log 2 + e^{-b} + \frac13
e^{-3b}+o\left(e^{-3b}\right)\right).
\label{btoa:eq}
\end{equation}
We can glue together four identical pentagons $Q$ and obtain a hyperbolic 
right-angled octagon $\widetilde Q$ as shown in Figure~\ref{pent:fig}(a). The octagon is
uniquely determined up to an isometry by $b$, which can vary freely in
$(0,+\infty)$. 

For visualisation purposes, consider the octagon $\widetilde Q$ as an
ordinary right-angled octagon on $\mathbb R^2$ plane, with four quarters of a
circle as alternating sides and other four sides parallel to coordinate axis. We
assign labels $\medtriangledown$ and $\medtriangleup$ to two sides parallel to the
horizontal axis and  $\medtriangleright$ and $\medtriangleleft$ to the sides parallel
to the vertical axis.
Translating the octagon along vertical and horizontal axes we obtain a tessellation
of a holed plane~$\Omega$, where the holes are Euclidean disks made of four quarters of
the boundary circles glued together, as shown in Figure~\ref{pent:fig}(b). 
The holed plane~$\Omega$ is a universal cover of a symmetric one-holed torus,
and carries the hyperbolic structure of~$\widetilde Q$. 
Evidently there is a natural action of the group $\Gamma = \mathbb Z \times
\mathbb Z$ on~$\Omega$ by isometries and copies~$\widetilde Q_{i,j}$ of~$\widetilde Q$ are
fundamental domains. Two  dashed lines in~$\widetilde Q$ give a pair of
shortest closed geodesics in~$\Omega/\Gamma$, which are generators of the
fundamental group.

We can use the local isometry between the plane and the torus with hyperbolic metric in order to
estimate the lengths of closed geodesics. 

It is known that every closed oriented geodesic $\gamma$ on $\Omega/\Gamma$ is freely homotopic to a
periodic word $\ldots \medtriangleup^{n_1}\medtriangleright^{k_1}
\medtriangleup^{n_2}\medtriangleright^{k_2} \ldots
\medtriangleup^{n_t}\medtriangleright^{k_t}\ldots$ of period 
$\omega(\gamma)\eqdef|n_1|+\ldots+|n_t|+|k_1|+\ldots+|k_t|$
where $n_j$ and $k_j$ are integers and $\medtriangleleft \! = \!
\medtriangleright^{-1}$, $\medtriangleup\! = \!\medtriangledown^{-1}$. We denote this geodesic by 
$$
\gamma_{( \medtriangleup^{n_1}\medtriangleright^{k_1}
\medtriangleup^{n_2}\medtriangleright^{k_2}\ldots
\medtriangleup^{n_t}\medtriangleright^{k_t})},
$$
and we call the periodic sequence $( \medtriangleup^{n_1}\medtriangleright^{k_1}
\medtriangleup^{n_2}\medtriangleright^{k_2}\ldots
\medtriangleup^{n_t}\medtriangleright^{k_t})$ the cutting sequence associated to
the closed geodesics.
A good exposition on cutting sequences associated to
closed geodesics on a one-holed torus can be found in~\cite{S85}.

The homotopy is unique up to conjugation by fundamental group. In
particular, the two shortest geodesics correspond to the ``constant'' sequences of period one
$(\medtriangleright)$ and $(\medtriangledown)$.  

The homotopy defines a bijection between closed geodesics and periodic two-sided infinite sequences
$\{\sigma_k\}_{-\infty}^{\infty}$ in the alphabet
$\Sigma=\{\medtriangleup, \medtriangledown,
\medtriangleright, \medtriangleleft \}$ which satisfy an additional condition
that $\sigma_k \ne \sigma_{k+1}^{-1}$
for any $k \in \mathbb Z$. We define a transition matrix
$$
\mathbb A = \left(\begin{matrix}
  1 & 0 & 1 & 1 \\
  0 & 1 & 1 & 1 \\
  1 & 1 & 1 & 0 \\
  1 & 1 & 0 & 1  \\ 
\end{matrix} \right).
$$
Let $\Sigma^{\mathbb A}$ be the set of words which are freely homotopic to
geodesics on $\Omega/\Gamma$, and periodic words in $\Sigma^{\mathbb A}$
correspond to closed geodesics. The shift on $\Sigma^{\mathbb A}$ corresponds to
the action of $\Gamma$ on closed geodesics.

It is easier to do the calculations using the upper half model $H = \{z \in \mathbb
C \mid \Im z > 0 \}$ of the hyperbolic plane and a subgroup of $PSL(2,\mathbb
R)$ acting on $H$. Namely, consider matrices
$$
B = 
\left(\begin{matrix}
    \phantom{-}\cosh(a) & -\sinh(a) \\
    -\sinh(a) & \phantom{-}\cosh(a)
\end{matrix} \right)
\qquad \qquad 
C = \left(
\begin{matrix}
  e^{a} & 0 \\
  0 & e^{-a} 
\end{matrix} \right)
$$
Then the subgroup $\langle B,C\rangle \subset PSL(2,\mathbb R)$ is a deck group
of the universal cover $H \to \Omega/\Gamma$ and generators $B$ and $C$
correspond to the generators $\gamma_{(\medtriangledown)}$ and
$\gamma_{(\medtriangleright)}$ of the fundamental group of~$\Omega/\Gamma$.
Moreover, the hyperbolic length of the geodesic corresponding to the cutting sequence  
$\ldots \medtriangleup^{n_1}\medtriangleright^{k_1}
\medtriangleup^{n_2}\medtriangleright^{k_2} \ldots
\medtriangleup^{n_t}\medtriangleright^{k_t}\ldots$ of period 
$\omega(\gamma)\eqdef|n_1|+\ldots+|n_t|+|k_1|+\ldots+|k_t|$, (where $k_j\ne0$,
$n_j\ne0$ for $j=1,\ldots,t$) is given by 
\begin{equation}
\ell\left(\gamma_{(\medtriangleup^{n_1}\medtriangleright^{k_1}
\medtriangleup^{n_2}\medtriangleright^{k_2} \ldots
\medtriangleup^{n_t}\medtriangleright^{k_t})} \right) = 
2 \acosh\left(\frac12\left|\tr\left(B^{n_1}C^{k_1}B^{n_2}C^{k_2} \ldots B^{n_t}
C^{k_t}\right)\right|\right).
\label{length:eq}
\end{equation}

\subsubsection{Approximating determinant}
\begin{notation}
  Given a contiguous subsequence $\sigma_1,\ldots,\sigma_k$ of a sequence $\sigma \in \Sigma^{\mathbb A}$ 
  we denote by $\gamma_{\sigma_1,\ldots,\sigma_k}$  a geodesic whose
  cutting sequence contains the contiguous subsequence $\sigma_1,\ldots,\sigma_k$.
  We denote by $\gamma_{[\sigma_1,\ldots,\sigma_k]}$ a segment of a geodesic whose
  cutting sequence contains a contiguous subsequence $\sigma_1,\ldots,\sigma_k$ with end
  points at the midpoints of the segments enclosed between intersections with the
  sides $\sigma_1,\sigma_2$ and $\sigma_{k-1},\sigma_{k}$.
  We denote by $\underline{\gamma_{[\sigma_1,\ldots,\sigma_k]}}$ the segment of the shortest
  of all closed geodesics whose cutting sequence contains the contiguous subsequence
  $\sigma_1,\ldots,\sigma_k$ with end points at the midpoints of segments enclosed
  between intersections with the sides labelled $\sigma_1$, $\sigma_2$ and
  $\sigma_{k-1}$, $\sigma_k$, respectively (see Figure~\ref{geod1:fig}). 
\end{notation}


Let $\sigma^1, \ldots, \sigma^N$ be all contiguous subsequences of sequences $\sigma \in
\Sigma^{\mathbb A}$ of length $n$. Let us consider an $N\times N$ transition
matrix given by 
$$
\mathbb A^n_{i,j} = 
\begin{cases}
  1, & \mbox{ if } \sigma^i_{k+1} = \sigma^j_k; \mbox{ for all } k = 1,2,\ldots,n-1. \\
  0, & \mbox{ otherwise.} 
\end{cases}
$$
We now define a one-parameter family of $N\times N$ matrices $A(s)$ which
elements depend on the length of geodesic segments determined by
$\sigma^i$ and $\sigma^j$: 
\begin{equation}
  A \colon \mathbb C \to \Mat(N\times N) \qquad A_{i,j}(s) = \mathbb A^n_{i,j} \cdot
  \exp(-s\cdot\ell(\gamma_{\underline{[\sigma^i_1\sigma^i_2\ldots\sigma^i_n\sigma^j_n]}})).
  \label{As-def:eq}
\end{equation}
Note that $A$ depends on $n$, but we omit this in the notation.
We will be using the following Lemma from~\cite{PV18} to find the approximate
location of the zeros.
\begin{lemma}
  Using the notation introduced above, we have the following representation for
  the Selberg zeta function
  \begin{equation}
    Z_{ \mathring{\mathbb T}(a) }(s)= \prod_{n=0}^\infty\prod_{\stackrel{\gamma = \mbox{\scriptsize
  primitive}}{\mbox {\scriptsize closed   geodesic}} } \left(1-
e^{-(s+n)\ell(\gamma)}\right) = \lim_{n\to \infty}\det(Id_n - A(s)),
  \end{equation}
  where $Id_n \in \mathrm{Mat}(n,n)$ is the identity matrix.
\end{lemma}
Our first Lemma gives approximations to lengths of geodesic segments.
\begin{lemma}
  \begin{align}
  \ell(\underline{\gamma_{[\medtriangledown\medtriangledown\medtriangleright]}})
  &= 2a - \log \sqrt 2 +  o(e^{-3a})& & \mbox{ see Figure~\ref{geod2:fig}(a)};\\
  \ell(\underline{\gamma_{[\medtriangleup\medtriangleright\medtriangleup]}}) &=
  2a - \log  2 + 2 e^{-2a} + o(e^{-3a}) & & \mbox{ see
  Figure~\ref{geod2:fig}(b)}; \\ 
  \ell(\underline{\gamma_{[\medtriangleup\medtriangleright\medtriangledown]}})
  &= 2a - \log  2 - 2e^{-2a} + o(e^{-3 a})& & \mbox{ see Figure~\ref{geod1:fig}(b)}.
  \end{align}
\end{lemma}
\begin{proof}
 The result follows by straightforward calculation by applying formula~\eqref{length:eq}
  together with~\eqref{btoa:eq} 
  to
  $\gamma_{(\medtriangleup\medtriangleright\medtriangledown\medtriangleleft)}$,
$\gamma_{(\medtriangledown\medtriangledown\medtriangleright\medtriangleright)}$,
and
$\gamma_{(\medtriangleup\medtriangleright\medtriangleup\medtriangleright)}$,
respectively. More precisely, by definition we have that
\begin{multline}
\ell(\underline{\gamma_{[\medtriangleup\medtriangleright\medtriangledown]}})
= \frac14
\ell(\gamma_{(\medtriangleup\medtriangleright\medtriangledown\medtriangleleft)})
= 
\frac12 \acosh \left(\frac12 \left|\tr (B^{-1}CBC^{-1})\right|\right) = \\ \frac12 \acosh\left(-
\frac{e^{4a}+e^{-4a}}{8} + \frac{e^{2a} + e^{-2a}}{2}   + \frac14 \right)
  = 2a -  \log 2 - 2 e^{-2a} - 5 e^{-4a} + o\left(e^{-4a}\right).
\end{multline}
Similarly,
\begin{multline}
 \ell(\underline{\gamma_{[\medtriangledown\medtriangledown\medtriangleright]}})
  = \frac14
  \ell(\gamma_{(\medtriangledown\medtriangleright\medtriangleright\medtriangledown)})
= \frac12 \acosh \left(\frac12\left| \tr (B B CC)\right| \right) = \\ 
\frac12 \acosh \left( \cosh^2 2a \right) 
  = 2a -  \log \sqrt 2 + 3 e^{-4a} + o\left(e^{-6a}\right);
\end{multline}
And in the last case 
\begin{multline}
  \ell(\underline{\gamma_{[\medtriangleup\medtriangleright\medtriangleup]}})
  = \frac14 \ell(\gamma_{(\medtriangleup\medtriangleright\medtriangleup\medtriangleright)})
= \frac12 \acosh \left(\frac12 \left|\tr (B^{-1}CB^{-1}C)\right|\right) = \\ \frac12 \acosh\left(
\frac{e^{4a}+e^{-4a}}{8} + \frac{e^{2a} + e^{-2a}}{2}   - \frac14 \right)
  = 2a -  \log 2 + 2 e^{-2a} - 5 e^{-4a} + o\left(e^{-4a}\right).
\end{multline} 

\begin{figure}
  \begin{tabular}{ccc}
\includegraphics{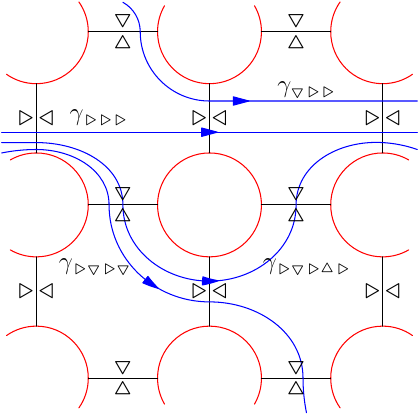} & \qquad \qquad & \includegraphics{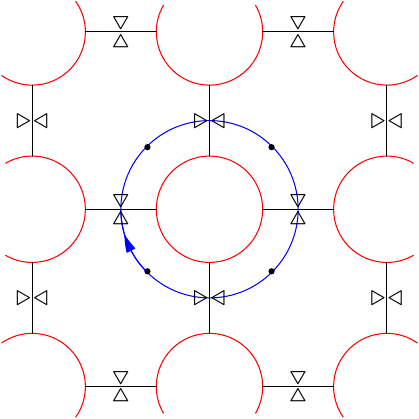}  \\
(a) & \ & (b) 
\end{tabular}
\caption{(a) Segments of geodesics lifted to the $\mathbb Z \times \mathbb Z$
holed plane, marked by contiguous subsequences of their cutting sequences, according to
visible intersections.
(b) The shortest closed geodesic among
$\gamma_{\medtriangleup\medtriangleright\medtriangledown}$
corresponding to periodic sequence
$(\medtriangleup\medtriangleright\medtriangledown\medtriangleleft)$
of period~$4$. The four marked points divide the geodesic into 4 equal segments: 
$\underline{\gamma_{[\medtriangleup\medtriangleright\medtriangledown]}}$,
$\underline{\gamma_{[\medtriangleright\medtriangledown\medtriangleleft]}}$,
$\underline{\gamma_{[\medtriangledown\medtriangleleft\medtriangleup]}}$,
and $\underline{\gamma_{[\medtriangleleft\medtriangleup\medtriangleright]}}$.
}
\label{geod1:fig}
\end{figure}

\begin{figure}
  \begin{tabular}{ccc}
    \includegraphics{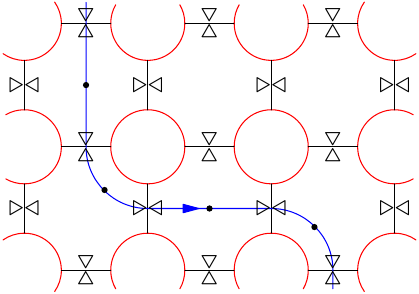} & \qquad \qquad & \includegraphics{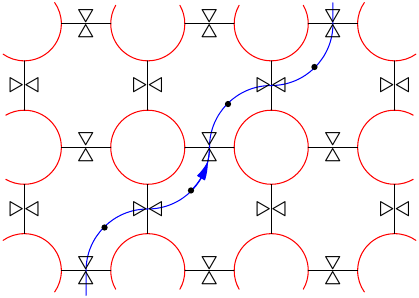} \\
    (a) & \ & (b)
  \end{tabular}
\caption{(a) The shortest closed geodesic among
$\gamma_{\medtriangledown\medtriangleright\medtriangleright}$,
corresponding to the periodic sequence
$(\medtriangledown\medtriangledown\medtriangleright\medtriangleright)$
of period $4$. The four marked points divide the geodesic into 4 equal segments: 
$\underline{\gamma_{[\medtriangledown\medtriangledown\medtriangleright]}}$,
$\underline{\gamma_{[\medtriangledown\medtriangleright\medtriangleright]}}$,
$\underline{\gamma_{[\medtriangleright\medtriangleright\medtriangledown]}}$,
and
$\underline{\gamma_{[\medtriangleright\medtriangledown\medtriangledown]}}$.
(b) The shortest closed geodesic among 
$\gamma_{\medtriangleup\medtriangleright\medtriangleup}$, 
corresponding to the periodic sequence
$(\medtriangleup\medtriangleright\medtriangleup
\medtriangleright)$ of period $4$.
 The four marked points divide the geodesic into 4 equal segments: 
$\underline{\gamma_{[\medtriangleup\medtriangleright\medtriangleup]}}$,
$\underline{\gamma_{[\medtriangleright\medtriangleup\medtriangleright]}}$,
$\underline{\gamma_{[\medtriangleup\medtriangleright\medtriangleup]}}$,
and
$\underline{\gamma_{[\medtriangleright\medtriangleup\medtriangleright]}}$.}
\label{geod2:fig}
\end{figure}
\end{proof}
In particular, we have the following corollary.
\begin{cor}
  \label{cor:length}
  $$
  \ell(\underline{\gamma_{[\sigma_1\sigma_2\sigma_3]}}) = \begin{cases}
    2a, &\mbox{ if } \sigma_1=\sigma_2=\sigma_3;\\
    2a - \log \sqrt 2 + o(e^{-3 a}), &\mbox{ if }  \sigma_1 = \sigma_2
    \ne \sigma_3; \\
    2a - \log 2 - 2e^{-2a} + o(e^{-3 a}), & \mbox{ if } \sigma_1 = \sigma_3
    \ne \sigma_2; \\
    2a - \log 2 + 2e^{-2a} + o(e^{-3 a}), &\mbox{ otherwise. }  \\
  \end{cases}
  $$
\end{cor}
\begin{remark}
  This explains why the case $a = k \log 2$ is different. In particular, we
  see that this choice makes the length of short geodesics rationally dependent. 
\end{remark}
We now apply Corollary~\ref{cor:length} to compute the matrix $A(s)$ for $n=2$. There are $12$
subsequences of length $2$ of the sequences from $\Sigma^{\mathbb A}$. We can
enumerate them as follows $\sigma^1= \medtriangleup\medtriangleup$,
$\sigma^2=\medtriangleup\medtriangleright$,
$\sigma^3=\medtriangleup\medtriangleleft$,
$\sigma^4= \medtriangleright\medtriangleright$, 
$\sigma^5=\medtriangleright\medtriangleup$,
$\sigma^6=\medtriangleright\medtriangledown$,
$\sigma^7= \medtriangledown\medtriangledown$, 
$\sigma^8=\medtriangledown\medtriangleright$,
$\sigma^9=\medtriangledown\medtriangleleft$,
$\sigma^{10}= \medtriangleleft\medtriangleleft$, 
$\sigma^{11}=\medtriangleleft\medtriangleup$,
$\sigma^{12}=\medtriangleleft\medtriangledown$. 
Then using definition~\eqref{As-def:eq} we may compute, for example 
\begin{align*}
  A_{1,1}(s) & = \exp(-s \cdot
  \ell(\underline{\gamma_{[\medtriangleup\medtriangleup\medtriangleup]}}))=\exp(-2as); \\
  A_{2,4}(s) & = \exp(-s \cdot 
  \ell(\underline{\gamma_{[\medtriangleup\medtriangleright\medtriangleright]}}))=\exp(-2as)\cdot\sqrt{2^s}
  \cdot\exp(-s\cdot o(e^{-3a}));  \\
  A_{2,5}(s) & = \exp(-s \cdot 
  \ell(\underline{\gamma_{[\medtriangleup\medtriangleright\medtriangleup]}}))=\exp(-2as)\cdot2^s\cdot\exp\left(2e^{-2a} s\right)  \cdot\exp(-s\cdot o(e^{-3a})); \\ 
  A_{2,6}(s) & = \exp(-s \cdot
  \ell(\underline{\gamma_{[\medtriangleup\medtriangleright\medtriangledown]}}))=\exp(-2as)
  \cdot{2^s} \cdot \exp\left(- 2  e^{-2a} s\right) \cdot\exp(-s\cdot o(e^{-3a})); 
\end{align*}
the other elements are similar. Observe that for values of $s$ within the critical strip,
$0<\Re s < 0.2$, we have that $|\exp(-s\cdot o(e^{-3a}))-1| \le 2 e^{-3a}$ is
small for $a$ sufficiently large.  Introducing a short-hand
notation 
\begin{align}
p_1(s) &= \sqrt{2^s}; \label{p1:eq}\\
p_2(a,s) &= 2^{s+1} \cosh(2s e^{-2a}); \label{p2:eq}\\
p_3(a,s) & = 2^{s+1} \sinh(2s e^{-2a}); \label{p3:eq}
\end{align}
we may consider a matrix 
\setcounter{MaxMatrixCols}{20}
$$
P(s) = 
\begin{pmatrix} 
     1 & p_1 & p_1 & 0 & 0 & 0 & 0 & 0 & 0 & 0 & 0 & 0  \\
     0 & 0 & 0 & p_1 & \frac{p_2-p_3}{2} & \frac{p_2+p_3}{2} & 0 & 0 & 0 & 0 & 0 & 0  \\
     0 & 0 & 0 & 0 & 0 & 0 & 0 & 0 & 0 & p_1 & \frac{p_2-p_3}{2} &
     \frac{p_2+p_3}{2}  \\
     0 & 0 & 0 & 1 & p_1 & p_1 & 0 & 0 & 0 & 0 & 0 & 0  \\
     p_1 & \frac{p_2-p_3}{2} & \frac{p_2+p_3}{2} & 0 & 0 & 0 & 0 & 0 & 0 & 0 & 0 & 0  \\
     0 & 0 & 0 & 0 & 0 & 0 & p_1 & \frac{p_2-p_3}{2} & \frac{p_2+p_3}{2} & 0 & 0 & 0  \\
     0 & 0 & 0 & 0 & 0 & 0 & 1 & p_1 & p_1 & 0 & 0 & 0  \\
     0 & 0 & 0 & p_1 & \frac{p_2+p_3}{2} & \frac{p_2-p_3}{2} & 0 & 0 & 0 & 0 & 0 & 0  \\
     0 & 0 & 0 & 0 & 0 & 0 & 0 & 0 & 0 & p_1 & \frac{p_2+p_3}{2} &
     \frac{p_2-p_3}{2}  \\
     0 & 0 & 0 & 0 & 0 & 0 & 0 & 0 & 0 & 1 & p_1 & p_1  \\
     p_1 & \frac{p_2+p_3}{2} & \frac{p_2-p_3}{2} & 0 & 0 & 0 & 0 & 0 & 0 & 0 & 0 & 0  \\
     0 & 0 & 0 & 0 & 0 & 0 & p_1 & \frac{p_2+p_3}{2} & \frac{p_2-p_3}{2} & 0 & 0 & 0  \\
\end{pmatrix}.
$$
Then the matrix $A(s)$ can be considered as a small perturbation of $\exp(-2as)P(s)$.
By Lemma~2.8 from~\cite{PV18}, the determinant $\det(I- A(s))$ approximates the
Selberg zeta function. 
Since the matrices $\exp(-2as)P(s)$ and $A(s)$ are close, the determinant
$\det(I - \exp(-2as)P(s))$ is an approximation to the zeta function, too. 
We see that the function $\det(I - \exp(-2as)P(s))$ is an exponential sum in $s$, and
therefore an almost periodic function with modul $\langle a, \log 2, 2
e^{-2a}\rangle$ if $a \not\in \log 2 \mathbb Z$ and $\langle \log 2, 
4^{-k}\rangle$ otherwise. It follows from general theory of almost periodic functions
that its zeros form a point-periodic set in the sense of
Krein--Levin. The simplicity of the matrix $P$ allows us to get more information
on exact location of the zeros of the determinant. 
Evidently, the matrix $\exp(-2as)P(s)$ has an eigenvalue $1$ if and only if $\exp(2as)$
is an eigenvalue of $P(s)$. The eigenvalues of the matrix $P(s)$ have a closed
form. 

\begin{align}
  \lambda_{1,2} &= \pm p_3 \label{lam12:eq}\\
\lambda_{3,4} &=  \frac{1}{2}\left(1-p_2 \pm
\sqrt{(p_2+1)^2 -8p_1^2 }\right) \label{lam34:eq}\\
 \lambda_{5,6} &=
 \frac12\left(1+p_2\pm\sqrt{(p_2-1)^2 +8 p_1^2}\right).\label{lam56:eq}
\end{align}
Introducing a shorthand notation 
$\tilde p(p_1,p_2,p_3)= 
9 p_3 \cdot (3 p_1^2 - p_2) - 1$ \label{p-tilde:def}
we write the remaining eigenvalues, each of multiplicity~$2$, as follows:
\begin{align}
  \lambda_{7,8} &= \frac13\left(1 +
  \frac{3p_3 p_2 -1}{\sqrt[3]{\sqrt{\tilde p^2-(3 p_3 p_2-1)^3}-\tilde p}} 
  +\sqrt[3]{\sqrt{\tilde p^2-(3p_3p_2-1)^3}- \tilde p}  \right); \label{lam78:eq}\\ 
  \lambda_{9,10} &= \frac{1}{3} -  \frac{i \sqrt3+1}{6}
  \cdot \frac{ 3p_3p_2 -   1 }{ \sqrt[3]{\sqrt{
  \tilde p^2 - (3p_3 p_2 - 1)^3}   - \tilde p}}  
  +  \frac{ i \sqrt3-1}{6}\cdot \sqrt[3]{ \sqrt{ \tilde p^2 - (3p_3p_2 - 1)^3}
  - \tilde p}; \label{lam910:eq} \\
  \lambda_{11,12} &= \frac{1}{3} +  \frac{ i \sqrt3-1}{6}
  \cdot \frac{ 3p_3p_2 -  1 }{ \sqrt[3]{\sqrt{
  \tilde p^2 - (3p_3p_2 - 1)^3}   - \tilde p}}  
  -  \frac{ i \sqrt3+1}{6}\cdot \sqrt[3]{ \sqrt{ \tilde p^2 - (3p_3p_2 - 1)^3} -
  \tilde p}. \label{lam1112:eq}
\end{align}
Summing up, we deduce the following. 
\begin{proposition}
\label{prop:zeros-loc}
Any {\bf small} zero of the zeta function $Z_{\mathring{\mathbb{T}}(a)}$ is close to a
solution of one of the twelve equations 
\begin{equation}
  \label{eq:zeros-loc}
\exp(2as)=\lambda_j(p_1(s),p_2(s),p_3(s)), \qquad j=1,\ldots,12,
\end{equation}
where $\lambda_k$,
$k=1,\ldots,12$ are given by~\eqref{lam12:eq}--\eqref{lam1112:eq}.
\end{proposition}
We omit the proof here. Figure~\ref{twoSets:fig} shows the small zeros of the zeta
function along with the zeros of the determinant.

In the next section we discuss properties of solutions of
equations~\eqref{eq:zeros-loc}, or in other words, zeros of the
determinant~$\det(I-\exp(-2as)P(s))$. 

\begin{figure}
\caption{Zeros of the determinant $I - \exp(-2as)P(s)$ (circles) and zeros of
the zeta function (dots) in a small part of the critical strip $1300 <\Im z <
1340$. In the cases we consider $\exp(\frac{3a}2)\approx 1400$.}
\begin{center}
\includegraphics{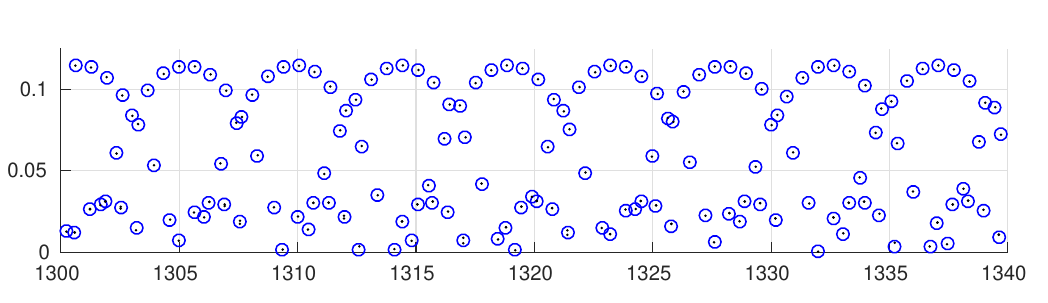}

(a) The case of $\mathring{\mathbb T}(10)$. 

\includegraphics{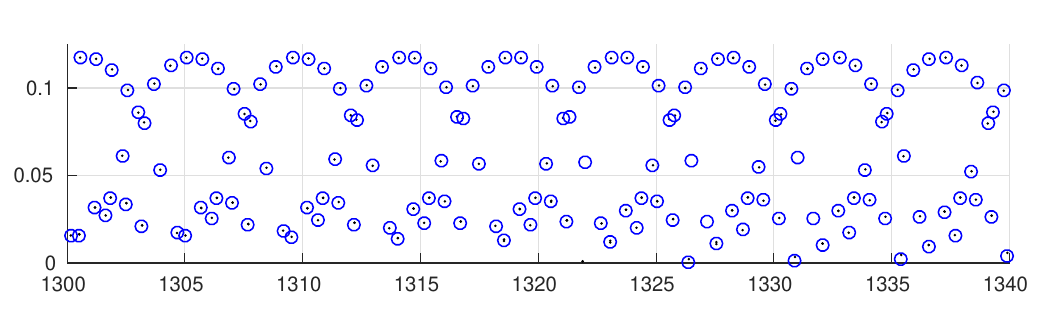}

(b) The case of $\mathring{\mathbb T}(14 \log 2+0.05)$. 

\includegraphics{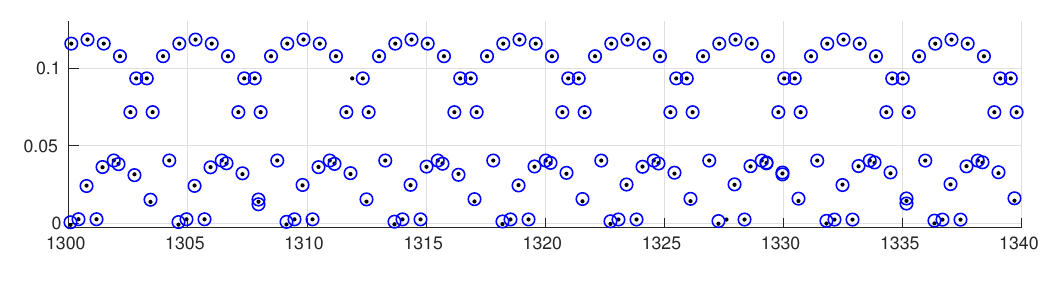}

(c) The case of $\mathring{\mathbb T}(14 \log 2)$. 

\end{center}
\label{twoSets:fig}
\end{figure}

\subsubsection{Zeros of the determinant}
\label{s:args}
\paragraph{Solving $\lambda_{1,2}(s)=\exp{2as}$.}
The first two eigenvalues have a simple form as functions of $s$ and
equation~\eqref{eq:zeros-loc} gives us two equations:
\begin{equation}
  \exp(2as) = \pm 2^{s+1}\sinh(2s e^{-2a}). 
  \label{lam12r:eq}
\end{equation}
We may write $s=\sigma+it$ and use the equality between squares of the absolute
values 
$$
|\exp\left( \left(2a-\log
2\right)s\right)|^2=\left|\sinh(2se^{-2a})\right|^2
$$
to obtain 
$$
\exp( (2a-\log 2)2 \sigma ) = 
\exp(4 \sigma e^{-2a})+\exp(-4 \sigma e^{-2a}) - 2 \cos(4t e^{-2a}), 
$$
which implies
\begin{equation}
t=\frac{e^{2a}}{4}\left( \pm \arccos \left( - \frac12  \exp( (2a - \log 2)2
\sigma ) +  \cosh (4 \sigma
e^{-2a}) )\right) + 2\pi k \right),
\label{tlam12sol:eq}
\end{equation}
provided 
$$
\left| \frac12  \exp( (2a - \log 2)2 \sigma ) -  \cosh (4 \sigma
e^{-2a}) \right| \le 1.
$$
Since for small $\sigma$ we have that $\cosh (4 \sigma
e^{-2a}) = 1 + 4 \sigma^2e^{-4a} + o\left(e^{-6a}\right)$, the above condition
is valid provided 
$ \sigma < \frac{\log 2}{2a - \log 2}$. However, for positive real part $\sigma>0$ we have
that imaginary part $t>\frac{1}4 e^{2a}$, which is outside of the range of small
zeros. Hence we have established the following.
\begin{lemma}
The first two eigenvalues~\eqref{lam12:eq} don't give any information on
location of small zeros in positive half-plane. 
\end{lemma} 
\begin{figure}
  \centerline{
  \includegraphics{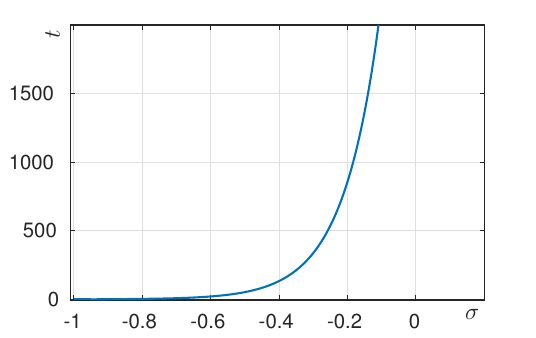}}
  \caption{A plot of imaginary part $\Im(s) = t$ as a function of the real part
  $\Re(s)=\sigma$ as defined by~\eqref{tlam12sol:eq}. The zeros of the
  determinant $\det(I-A(s))$ defined by $\lambda_{1,2}=\exp(2as)$ are outside of
  the domain of approximation. } 
  \label{lam12:fig}
\end{figure}

The equality~\eqref{lam12r:eq} also implies the
equality between arguments:
$$
\Arg\left( \exp\left( \left(2a-\log 2\right)s\right) \right)=
\Arg\left(\sinh(2se^{-2a})\right).
$$
We see that for $s=\sigma_0 + it$ the function
$\Arg\left(\sinh(2se^{-2a})\right)$ is monotone increasing and small, while
$\Arg\left( \exp\left( \left(2a-\log 2\right)s\right) \right)$ is changing
rapidly. We therefore expect that solutions of~\eqref{lam12r:eq} with
$|s|<e^{2a}$ belong to the 
curve given by~\eqref{tlam12sol:eq} and the difference between imaginary parts
of consecutive zeros is approximately $\frac{2\pi}{2a - \log 2}$.

\begin{figure}
  \caption{Families of curves $\mathcal T_1$ quickly oscillating with period
  $\frac{\pi}{a}$  and $\mathcal T_2$ slowly oscillating with period
  $\frac{2\pi}{\log 2}$ as defined by~\eqref{eq:T1def} and~\eqref{eq:T2def} and
  described in Lemma~\ref{lem:lam34} and zeros of the 12th Taylor polynomial 
  approximating the zeta function (stars). We see that
  the actual zeros occur very close to odd elements of the sequences of points
  of intersections.  
  Zeros on the imaginary axis correspond to solutions of $e^{2as} = \lambda_{7,8}(s)$. }
  \vspace*{-7mm}
\begin{center}
  \includegraphics{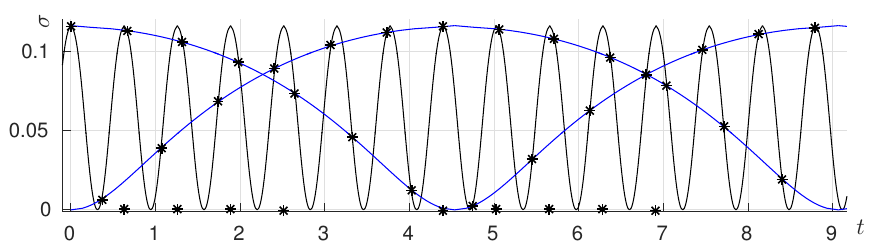} \\
  (a) Case $a=10$. A part of the zero set from Figure~\ref{fig:torus-plots}(a).
  \includegraphics{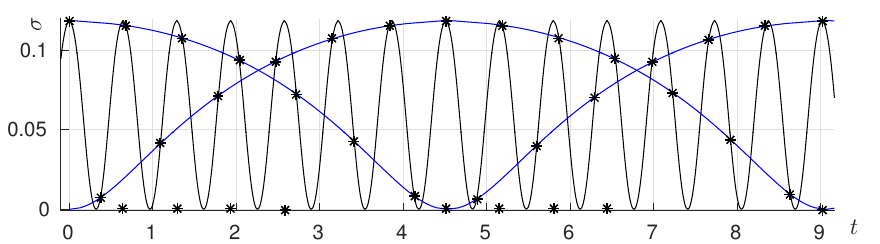} \\
  (b) Case $a=14\log2+0.05$. A part of the zero set from
  Figure~\ref{fig:torus-plots}(b).
  \includegraphics{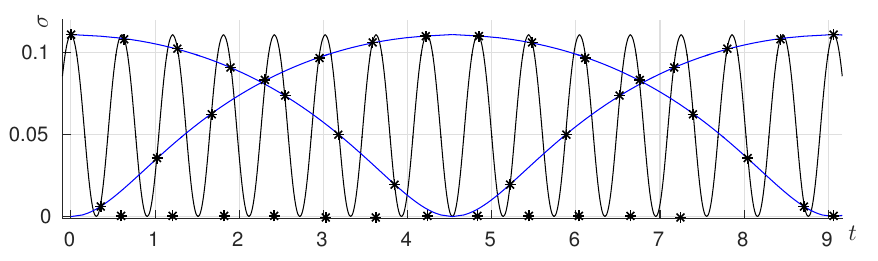} \\
  (c) Case $a=14 \log 2$. A part of the zero set from
  Figure~\ref{fig:torus-plots}(c). The real and the imaginary axis are swapped.
  The nearly straight lines would correspond to $\sigma=\const$. 
  \includegraphics{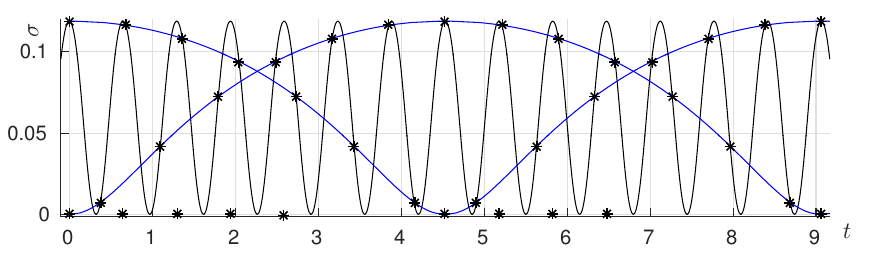} \\
  (d) Case $a=15 \log 2$. Together with the plot above they illustrate Corollary~\ref{cor:lam34}.
\end{center}
 \label{fig:Tfamily}
\end{figure}

We now proceed to analyse the equations coming from the next four eigenvalues.
\paragraph{Solving $\lambda_{3,4}(s)=\exp{2as}$ and
$\lambda_{5,6}(s)=\exp{2as}$.}
The equations~\eqref{eq:zeros-loc} read 
\begin{align}
  2 e^{2as} &= 1 - 2^{s+1} \cosh(2s e^{-2a}) \pm \sqrt{(2^{s+1}\cosh(2s e^{-2a})
  + 1)^2 - 2^{s+3}}; \label{lam34-2:eq} \\
  2 e^{2as} &= 1 + 2^{s+1} \cosh(2s e^{-2a} ) \pm \sqrt{(2^{s+1} \cosh(2s
  e^{-2a}) - 1)^2 + 2^{s+3}}; \label{lam56-2:eq}
\end{align}
which is equivalent to
\begin{align*}
e^{2as}\left(e^{2as}\cdot 2^{-s-1} - 2^{-s-1} + \cosh(2s e^{-2a})\right) =
\cosh(2s e^{-2a}) -1, \\
e^{2as}\left(e^{2as}\cdot 2^{-s-1} - 2^{-s-1} - \cosh(2s e^{-2a})\right) =
1-\cosh(2s e^{-2a}).  
\end{align*}
The determinant $\det(I - \exp(-2as)P(s))$ approximates the zeta function on a
part of the critical strip $0<\Re(s)<\delta$, $0<\Im(s)<e^{\frac{3a}2}$.
We have the following asymptotic expansion for the right hand side:
\begin{multline}
\left|\cosh(2s e^{-2a}) - 1 \right| \le \sum_{j=1}^\infty \left|4^j e^{-4aj}
s^{2j} \right| = \sum_{j=1}^\infty 4^j e^{-4aj} (\sigma^2 + t^2)^j \\ \le
4e^{-4a}(\sigma^2 + t^2) + \frac{16 e^{-8a} (\sigma^2 + t^2)^2}{1-4 (\sigma^2 +
t^2) e^{-4a}}. 
\end{multline}
This allows us to deduce that zeros of the zeta function with imaginary part
$\Im(s) = |t|\le e^{\frac{3a}2} \approx 1800$ are close to
solutions of the approximate equations 
\begin{align}
  e^{2as}\cdot 2^{-s-1} - 2^{-s-1} + 1  = 0, \label{lam34approx:eq} \\
  e^{2as}\cdot 2^{-s-1} - 2^{-s-1} - 1  = 0. \label{lam56approx:eq}
\end{align}
Evidently solutions of~\eqref{lam34approx:eq} and~\eqref{lam56approx:eq} should satisfy 
\begin{align}
  |e^{2as}|^2 & = |2^{s+1} - 1|^2 \mbox { and }  |2^{s+1}|^2 = |e^{2as}-1|^2;
  \mbox{ or } \\
  |e^{2as}|^2 &= |2^{s+1} + 1|^2 \mbox { and }  |2^{s+1}|^2 = |e^{2as}-1|^2. 
\end{align}
and therefore belong to the intersections $\mathcal T_1 \cap \mathcal T_2$ 
of the of curves given by 
\begin{align}
  \mathcal T_1 & \eqdef \left\{ \sigma + it \, \Bigl| \, \cos(2at) =
  \pm \frac{1 + e^{2a\sigma} - 4^{1+2\sigma}}{2 e^{2a\sigma}} \right\}
  \label{eq:T1def} \\
  \mathcal T_2  & \eqdef \left\{ \sigma + it \, \Bigl| \, \cos(t \log 2 ) = \pm \frac{4^{1+2\sigma} + 1 -
  e^{4a\sigma}}{4^{1+\sigma}} \right\}. \label{eq:T2def} 
\end{align}
We summarize our fundings in the following
\begin{lemma}
  \label{lem:lam34}
There exist a zero of the zeta function~$\zeta_{\mathring{\mathbb T}(a)}$ in 
$e^{-a}$-neighbourhood of every \textbf{odd} element of the following subsequences of the 
points of intersection with $\Im z_n < e^{\frac{3a}{2}}$ (see plots in Figure~\ref{fig:Tfamily})
\begin{align*}
 z_n &=\! \sigma_n + it_n \mbox{ where } \cos(2at_n) \!=\!
\frac{1 + e^{2a\sigma_n} - 4^{1+2\sigma_n}}{2 e^{2a\sigma_n}}, \ 
\cos(t_n \log 2 ) \!= \!\frac{e^{4a\sigma_n}-4^{1+2\sigma_n} - 1}{2^{2+\sigma_n}}
\mbox{ and } t_n<t_{n+1}; \\
z_n &=\! \sigma_n + it_n \mbox{ where } \cos(2at_n) \!=\!
\frac{4^{1+2\sigma_n}-1 - e^{2a\sigma_n}  }{2 e^{2a\sigma_n}}, \ 
\cos(t_n \log 2 ) \!= \!\frac{e^{4a\sigma_n}-4^{1+2\sigma_n} - 1}{2^{2+\sigma_n}}
\mbox{ and } t_n<t_{n+1}. 
\end{align*}
\end{lemma}
\begin{cor}
  \label{cor:lam34}
In the case $a \in \mathbb N \log 2 $ solutions to this system belong to the
straight lines $\Im z=\sigma=\const$; moreover, the intersection of the zero set
with any of this lines is a periodic set of period $\frac{2\pi}{\log 2}$. 
These lines correspond to seemingly straight lines we see in Figure~\ref{fig:torus-plots}(c). 
In the case $a \not\in \mathbb N\log 2 $ this no longer holds and we see a
random structure as shown in Figure~\ref{fig:torus-plots}(a). 
\end{cor}

\begin{figure}
  \caption{Several plots showing the curve $|\exp(2as)|=|\lambda_{7,8}(s)|$. (a)
  General plot of the curve within the rectangle $-0.01< \Re s < 0.01$, $0< \Im
  s < 200$; (b) A part of the curve near the real axis $-0.005< \Re s < 0.05$,
  $0 < \Im s < 20$, stars mark zeros of the zeta function, the oscillating curve in horizontal
  direction is $\mathcal T_2$; (c) The curve in the part of the critical strip
  $0<\Re s<\delta$, $310<\Im s <340$, oscillating around imaginary axis. Stars
  mark zeros of the zeta function. The second oscillating curve is $\mathcal
  T_2$. }
  \hspace*{-3mm}\begin{tabular}{ccc}
  \includegraphics[scale=0.75]{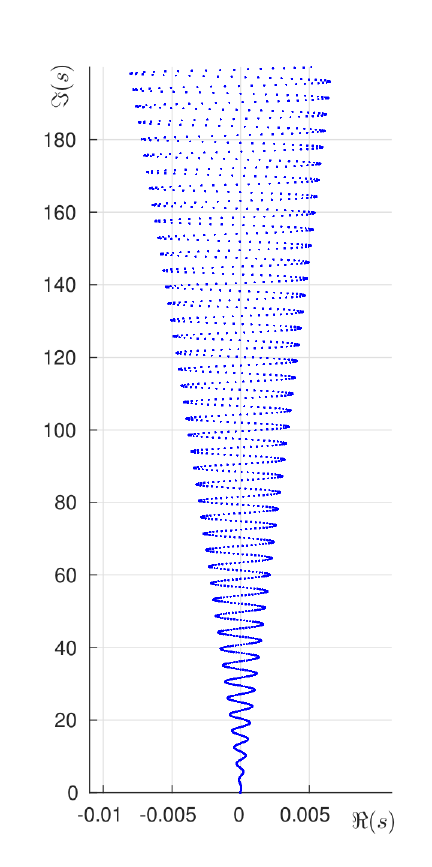} &
  \includegraphics[scale=0.75]{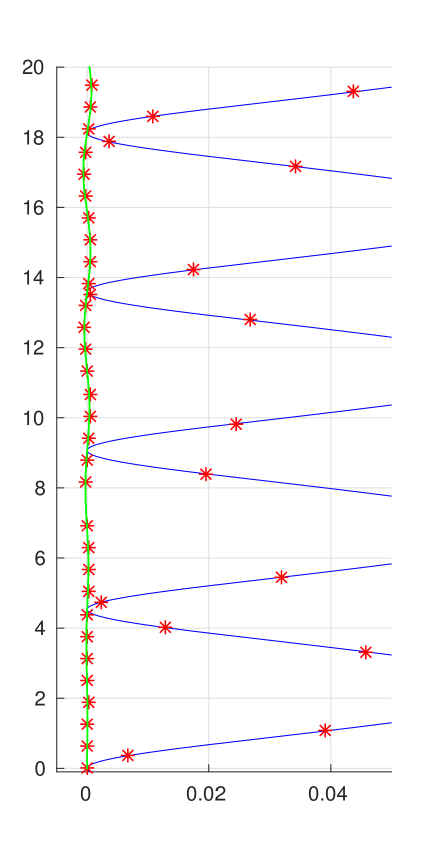} &
  \includegraphics[scale=0.75]{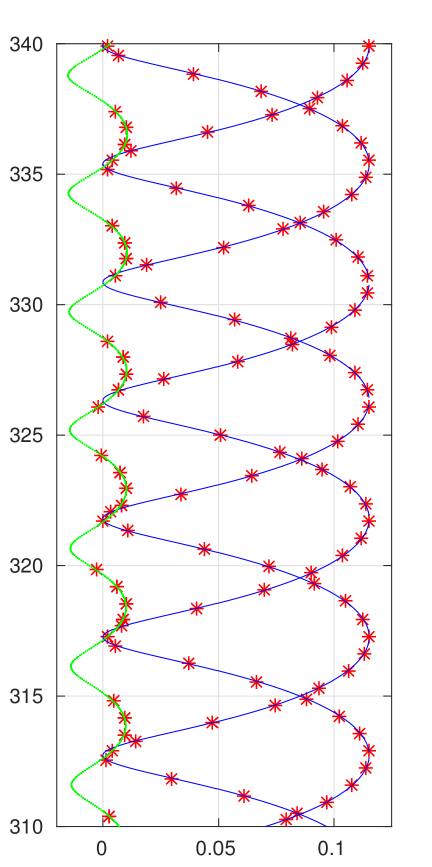} \\
  (a) & (b) & (c)
\end{tabular}
\label{fig:lam78}
\end{figure}

\paragraph{Analysing $\lambda_{7,8}(s)=\exp(2as)$}
Let us now consider the remaining eigenvalues $|\exp(2as)| = |\lambda_{k}(s)|$,
$7 \le k \le 12$. The explicit expressions~\eqref{lam78:eq},~\eqref{lam910:eq}
and~\eqref{lam1112:eq} can be used to compute the curves where the zeros are
located numerically. It turns out that equation~\eqref{eq:zeros-loc} have
solutions in the domain of approximation only for $k=7,8$. In
Figure~\ref{fig:lam78}(a) we see a part of the curve defined by 
\begin{equation}
\mathcal T_3=\left\{ s=\sigma+it \mid |\exp(2as)| = |\lambda_{k}(s)| \right\}
\label{eq:t3def}
\end{equation}
for $k=7,8$; Figures~\ref{fig:lam78}(b) and~\ref{fig:lam78}(c) show two
pieces of the curve corresponding to $0<\Im s<20$ and $310< \Im s< 340$ with
zeros of the zeta function located there. 

In an attempt to find a transversal family of curves that would help to describe
the location of the zeros more precisely, we would like to study the
polynomial 
$$
\prod_{j=7}^{12}(x-\lambda_j(s)) = \det (P(s) - x I) \cdot
\prod_{j=1}^{6}(x-\lambda_j(s))^{-1}.
$$
By a straightforward calculation we obtain 
\begin{equation}
\prod_{j=7}^{12}(x-\lambda_j(s)) = 
\left(x^3-x^2 + p_2 p_3 x + (2p_1^2-p_2) p_3  \right)^2.
\label{eq:cubicFactor}
\end{equation}
We would like to return to the original variable~$s$, i.e. to
reverse~\eqref{p1:eq}--\eqref{p3:eq}, taking into account that for
$\Re(s)\ll1$ we have that $\sinh(s) \approx i \sin(t)$ and $\cosh(s)
\approx \cos(t)$:
\begin{align}
  p_1(s) & \approx 2^{\frac s2}; \\
  p_2(s) & \approx 2^{s+1} \cos( 2  t e^{-2a}); \\
  p_3(s) & \approx  2^{s+1}  \sin (2 t e^{-2a} ) i.
\end{align}
We obtain an approximation 
\begin{equation*}
 \prod_{j=7}^{12}(x-\lambda_j(s))  \approx  
\left(x^3 - x^2 + i2^{2s+1} \left(  \sin(4 e^{-2a}t) (x-1) + 2\sin(2 e^{-2a}t)
\right)\right)^2.
\end{equation*}
The latter implies 
\begin{equation}
 2^{4\sigma+2} = \left| 2^{2s+1} \right|^2 =
\left|   \frac{(x^3 - x^2)}{ (\sin(4
  e^{-2a} t)(x-1)  + 2\sin(2 e^{-2a}t))  } \right|^2.
\label{eq:lam78}
\end{equation}
Substituting $x=\exp(2as)$ we solve the equation for $\cos(2at)$ and obtain the
equality
\begin{equation}
  \cos(2at)=\frac{e^{8a\sigma}(e^{4a\sigma}+1)-2^{4\sigma+2}\left(\sin^2(4e^{-2a}t)(e^{4a\sigma}+1)
  - 2\sin(4e^{-2a}t) + 4 \sin^2(2e^{-2a}t)\right) }{
   2 e^{2a\sigma} \left(2^{4\sigma+2} (\sin(4    e^{-2a}t) - \sin^2(4e^{-2a}t))
  +  e^{8a\sigma}\right)}.
\end{equation}
We see that although the right hand side depends on $t$, and as $t$ varies in any small interval
of length $c\ll \exp(a)$, say $(nc,(n+1)c)$ the dependence on $t$ is negligible,
so we may consider a partition into intervals of length $c$ and the curves defined by
\begin{multline}
  \mathcal T_4 = \biggl\{ \sigma+it \; \Bigl| \; \cos(2at) = \\  \frac{e^{8
  a\sigma} (e^{4a\sigma+1}) \!-\!
  2^{4\sigma+2}\left(\sin^2(4e^{-2a}nc)(e^{4a\sigma}+1)
  - 2\sin(4e^{-2a}nc) + 4 \sin^2(2e^{-2a}nc)\right) }{
   2 e^{2a\sigma} \left(2^{4\sigma+2} (\sin(4    e^{-2a}nc) - \sin^2(4e^{-2a}nc))
   +  e^{8a\sigma}\right)} \biggr\}
  \label{eq:t4def}
\end{multline}
on the intervals $nc<t<(n+1)c$, $n \in \mathbb Z$. 
\begin{remark} The dependence of the right
hand side on $t$ is reflected in increasing amplitude of oscillations of the
curve in Figure~\ref{fig:lam78}(a).  It is possible to make further
simplification of the right hand side of~\eqref{eq:t4def}, using first order approximations  
$\sin(x) \approx x$ and $\cos(x) \approx 1$ for small $x$. This would lead to 
\begin{equation}
 \mathcal T_4^\prime = \biggl\{ \sigma+it \, \Bigl| \, \cos(2at)= 
 \frac{e^{8a\sigma}(e^{4a\sigma}+1) - 2^{4\sigma+5} e^{-2a} t (2e^{-2a} t
 (e^{4a\sigma}+2)-1)}{2 e^{2a\sigma}(2^{4\sigma+4} e^{-2a} t (1-4 e^{-2a} t) +
 e^{8a\sigma} )} \biggr\}. 
\label{eq:t4dash}
\end{equation}
It is evident that the plot will be a quickly oscillating curve with period
$\frac{\pi}{a}$ and increasing amplitude of oscillations. 
\end{remark}
We summarize our discussion in this section as follows. 
\begin{proposition}
  Zeros of the zeta function in the critical strip with imaginary part $\Im s <
  e^a$ are $e^{-a}$-close either to the intersections $\mathcal T_1 \cap
  \mathcal T_2$, as described in Lemma~\ref{lem:lam34} or to the intersections
  $\mathcal T_3 \cap \mathcal T_4$.  
\end{proposition}

\end{document}